\documentclass[11pt]{article} 

\usepackage{amssymb}
\usepackage{amsmath}
\usepackage{geometry}
\usepackage{braket}
\usepackage{subfig}
\usepackage{url}
\usepackage{algorithm}
\usepackage{algpseudocode}
\usepackage{setspace}
\usepackage{xcolor}

\usepackage{graphicx} 
\usepackage{multirow}
\usepackage{amsthm}
\usepackage{graphicx}
\usepackage[linktoc=all]{hyperref}
\hypersetup{citecolor=linkcolor}

\newtheorem{theorem}{Theorem}
\newtheorem{lemma}[theorem]{Lemma}
\newtheorem{example}[theorem]{Example}
\newtheorem{corollary}[theorem]{Corollary}

\newtheorem{remark}[theorem]{Remark}
\newtheorem{definition}[theorem]{Definition}

\algrenewcommand\algorithmicrequire{\textbf{Input:}}
\algrenewcommand\algorithmicensure{\textbf{Output:}}

\title{Randomized low-rank Runge-Kutta methods\footnote{This work was supported by the SNSF research project \textit{Fast algorithms from low-rank updates}, grant number: 200020\_178806.}}

\author{Hei Yin Lam\thanks{\'Ecole Polytechnique F\'ed\'erale de Lausanne (EPFL), Institute of Mathematics, Switzerland. \texttt{hysan.lam@epfl.ch}, \texttt{daniel.kressner@epfl.ch}.} \and Gianluca Ceruti\thanks{Department of Mathematics, University of Innsbruck, Innsbruck, Austria. \texttt{gianluca.ceruti@uibk.ac.at}.} \and Daniel Kressner\footnotemark[2]}
\date{\today}

\begin{document}

\maketitle
\begin{abstract}
This work proposes and analyzes a new class of numerical integrators for computing low-rank approximations to solutions of matrix differential equation.
We combine an explicit Runge-Kutta method with repeated randomized low-rank approximation to keep the rank of the stages limited. The 
so-called generalized Nystr\"om method is particularly well suited for this purpose; it builds low-rank approximations from random sketches of the discretized dynamics.
In contrast, all existing dynamical low-rank approximation methods are deterministic and usually perform tangent space projections to limit rank growth. Using such tangential projections can result in larger error compared to approximating the dynamics directly. Moreover, sketching allows for increased flexibility and efficiency by choosing structured random matrices adapted to the
structure of the matrix differential equation.
Under suitable assumptions, we establish moment and tail bounds on the error of our randomized low-rank Runge-Kutta methods.
When combining the classical Runge-Kutta method with generalized Nystr\"om, we obtain a method called Rand RK4, which exhibits fourth-order convergence numerically -- up to the low-rank approximation error. For a modified variant of Rand RK4, we also establish fourth-order convergence theoretically. Numerical experiments for a range of examples from the literature
demonstrate that randomized low-rank Runge-Kutta methods compare favorably with two popular dynamical low-rank approximation methods, in terms of robustness and speed of convergence.
\end{abstract}

\section{Introduction}

In this work, we aim at approximating the solution $A(t)$ to large-scale matrix differential equations of the form
\begin{equation}
\label{ode}
    \dot{A}(t)=F(A(t)),\qquad A(0)=A_0\in \mathbb{R}^{m\times n}.
\end{equation}
In many situations of practical interest, an autonomous ordinary differential equation can be naturally viewed as such a matrix differential equation. Examples include applications
in physics~\cite{kieri2016discretized,trombettoni2001discrete}, uncertainty quantification~\cite{Babaee2017robust}, and machine learning~\cite{schotthofer2022low}. For large $m$, $n$, the solution of~\eqref{ode} becomes expensive; in fact, it may not even be possible to store the entire matrix $A(t)$ explicitly. To circumvent this limitation, model order reduction techniques can be employed. An increasingly popular approach is based on exploiting (approximate) low-rank structure of $A(t)$, which arises, for example, from smoothness properties of the underlying physical system.
In particular, dynamical low-rank approximation~\cite{Koch2007Dynamical} approximates $A(t)$ by evolving matrices $Y(t)$ on the manifold $\mathcal{M}_r$ of rank-$r$ matrices.
As only the rank-$r$ factors of $Y(t)$ need to be stored, this already reduces memory requirements significantly when $r\ll m,n$.
By the Dirac-Frenkel variational principle, the matrix $Y(t)$ is obtained by solving the differential equation
\begin{equation}
\label{dlra-ode}
    \dot{Y}(t)=P_r(Y(t))F(Y(t)),\quad Y(0)=Y_0\in \mathcal{M}_r,
\end{equation} where $P_r(Y(t))$ denotes the orthogonal projection onto $T_{Y(t)} {\mathcal M}_r$, the tangent space of ${\mathcal M}_r$ at $Y(t)$.
To also achieve a reduction of computational cost, one needs to exploit the low-rank structure of $Y(t)$ when integrating~\eqref{dlra-ode}. As shown in~\cite{Koch2007Dynamical}, this can be achieved by rewriting~\eqref{dlra-ode} as a system of differential equations for the rank-$r$ factors of $Y$. However, directly integrating this system with standard explicit time integration methods often leads to poor approximation results, unless (very) small time step sizes are used. This is caused by the additional stiffness introduced by small singular values of $Y(t)$. To address this issue, special integrators have been proposed that are robust to the presence of small singular values and allow for much larger step sizes. These methods include the projected splitting integrator \cite{Lubich2014projector}, projection methods \cite{kieri2019projection,Charous2023Dynamically} as well as Basis Update \& Galerkin (BUG) integrators \cite{ceruti2024robust,Ceruti2024,ceruti2022unconventional}. Under the assumption 
\begin{equation} \label{eq:assumption}
    \|F(Y)-P_r(Y)F(Y)\|_F\leq \tilde{\epsilon}, \text{ for all }Y\in \mathcal{M}_r\cap \{\text{suitable neighbourhood of $A(t)$}\} 
\end{equation} all these methods exhibit at least first-order convergence up to $\mathcal{O}(\tilde{\epsilon})$, both theoretically and numerically. The mid-point BUG~\cite{ceruti2024robust}, a variant of the parallel integrator~\cite{kusch2024second}, and the projected Runge–Kutta methods~\cite{kieri2019projection} are the only provable second-order integrators up to $\mathcal{O}(\tilde{\epsilon})$. Projected Runge–Kutta methods can also achieve higher order.

Assumption~\eqref{eq:assumption}, which says that $F(Y)$ is nearly contained in the tangent space, is arguably a strong assumption. For small $\tilde{\epsilon}>0$ it implies, at least for short times, that $A(t)$ can be well approximated by a rank-$r$ matrix, but the reverse is not true. In particular, it is possible that $A(t)$ can be well approximated by a rank-$r$ matrix even if~\eqref{eq:assumption} is not satisified with small $\tilde{\epsilon}$. According to~\cite{kieri2019projection}, this can occur, for instance, when the manifold $\mathcal M_r$ close to $A(t)$ has tiny, high-frequency wiggles, or when the range and co-range of $Y$ are contained in the orthogonal complement of the range and co-range of $F(Y)$, respectively~\cite{appelo2024robust}. Concrete examples are given in Section \ref{sect: Numerical Experiments}. When Assumption~\eqref{eq:assumption} is not satisifed with small $\tilde{\epsilon}$, the use of tangent space projections in numerical methods bears the danger of introducing unacceptably high errors.

In this work, we develop low-rank time integration methods for~\eqref{ode} that do not rely on~\eqref{eq:assumption} but only require $A(t)$ to admit accurate low-rank approximations. Our approach is based on the notion of projected integrators~\cite[Ch. IV.4]{Hairer2010Geometric}, which first perform a standard time integration step and then project back to the manifold. For the manifold $\mathcal M_r$, the efficiency of projected integrators is impaired by the occurrence of high-rank matrices, e.g., during the intermediate stages of a Runge-Kutta method. In~\cite{kieri2019projection}, this issue is addressed by repeatedly applying tangent space projection, which limits the rank to $2r$ at the expense of having to impose~\eqref{eq:assumption}.

In this work, we take a novel approach to avoid the high ranks encountered by projected integrators. Our approach is based on performing randomized low-rank approximation, which uses random sketches instead of tangent space projections. For a \emph{constant} matrix $B$, such randomized approaches have been studied intensively during the last 1-2 decades, including the popular randomized SVD~\cite{Halko2011Finding} and the (generalized) Nystr\"{o}m approximation \cite{nakatsukasa2020fast,Tropp2017Practical}.
The generalized Nystr\"{o}m approximation utilizes two sketches $B\Omega$ and $\Psi^TB$ to approximately capture the range and co-range of $B$, where $\Omega$ and $\Psi$ are random matrices with the number of columns chosen to be slightly larger than $r$. 

To the best of our knowledge, this is the first work that proposes and analyzes randomized low-rank approximation methods for time integration. The randomized low-rank Runge-Kutta (RK) methods proposed in this work combine explicit RK methods with randomized low-rank approximation. Our analysis applies to any randomized low-rank approximation satisfying the moment assumptions defined in Section~\ref{section: Randomized low-rank approximation}, but for simplicity we will restrict our algorithmic considerations to the generalized Nystr\"{o}m approximation. Assuming that the dynamics generated by $F$ preserve rank-$r$ matrices approximately, we derive a probabilistic result that establishes a convergence order (up to the level of rank-$r$ approximation error) based on the so-called stage order of the underlying RK method. This matches the order established in~\cite{kieri2019projection} for projected RK methods. However, unlike the results in~\cite{kieri2019projection}, our numerical experiments indicate that randomized low-rank RK methods actually achieve the usual convergence order of the RK method, which can be significantly higher. For the randomized low-rank RK method based on RK 4, we also establish order $4$ theoretically when allowing for modest intermediate rank increases in the stages. This compares favorably to order $2$ implied by the techniques from in~\cite{kieri2019projection}. 

The remainder of the paper is organized as follows. In Section \ref{Section: 2}, after providing some preliminaries, we propose and analyze an idealized projection method, which assumes that the exact flow of~\eqref{ode} is given. We show that applying randomized low-rank approximation causes, with high probability, little to no harm to time integration. In Section \ref{Section: 3}, we propose a practical method that uses an RK method to approximate the exact flow and provide error analysis based on the stage order. Furthermore, we prove that if we allow rank increase in the intermediate stages then the classical RK method combined with randomized low-rank approximation can still achieve convergence order 4, up to the level of the low-rank approximation error. Finally, in Section \ref{sect: Numerical Experiments}, we provide a range of numerical experiments that confirm the theoretical results and demonstrate the robust convergence of randomized RK methods.

\section{Preliminaries and an idealized randomized projection method}

\label{Section: 2}
In this section, we provide preliminaries on (randomized) low-rank approximation and introduce an idealized randomized  projection method that allows one to study the impact of randomization in an isolated fashion. In the following, $\|\cdot\|_F$ denotes the Frobenius norm of a (constant) matrix and $\|Y \|_{L_q}=(\mathbb{E}[\|Y\|_F^q])^\frac{1}{q}$ denotes the $L_q$ norm, for some $q\geq 1$, of a random matrix $Y$.

\subsection{Assumptions}
\label{section: assumptions}

Our analysis will be based on the following three assumptions on $F$. The first two assumptions are the same as in~\cite{ceruti2024robust, kieri2019projection, kusch2024second}, while the third one is a modification of the usual low-rank approximability assumption in dynamical low-rank approximation.

{\bf Assumption 1}. We assume that $F$ is Lipschitz continuous, that is, there is a Lipschitz constant $L>0$ such that
\begin{equation}
\label{eq:Lipschitz_assumption}
    \|F(X)-F(Y)\|_F\leq L\|X-Y\|_F\quad \text{for all} \quad X,Y\in \mathbb{R}^{m\times n}.\end{equation}
 By the Picard-Lindel\"{o}f theorem, this implies that the solution of~\eqref{ode} exists and it is unique for some finite time interval.

{\bf Assumption 2}. Let $\Phi^t_F$ denote the exact flow of $F$, that is, 
given the solution $A(t)$ of~\eqref{ode}, we have that $A(t)=\Phi^t_F(A_0)$. For a method of order $\tau$, we assume that the first $\tau$ derivatives
\begin{equation}
\label{eq:bounded_assumption}
    \frac{\mathrm{d}^{j+1}}{\mathrm{d}t^{j+1}}\Phi^t_F(Y), \quad j = 0,1,\ldots,\tau,
\end{equation} have uniformly bounded norm for all $Y\in \mathbb{R}^{m\times n}$. This assumption is needed when, e.g., performing local error analysis of higher-order methods~\cite[Chapter II.1]{Hairer1987Solving}.

{\bf Assumption 3}. To ensure low-rank approximability, we assume for every $h\leq h_0$
that
 \begin{equation}
\label{eq:low-rank_assumption}
    \|\Phi^h_F(Y)-[\![\Phi^h_F(Y)]\!]_r\|_{F}\leq C_Mh\epsilon \quad \text{for all} \quad Y\in \mathcal{M}_r,
\end{equation} where $C_M$ is a constant that depends on $L$ and $h_0$ only. Here and in the following, we use $[\![\cdot]\!]_r$ to denote a best rank-$r$ approximation of a matrix. 
Assumption~\eqref{eq:low-rank_assumption} replaces the assumption \eqref{eq:assumption}, usually made in the analysis of dynamical low-rank approximation~\cite{kieri2016discretized,kieri2019projection}; see Section~\ref{sec:lorapprox} below for a more detailed discussion.

In this paper, we assume that \eqref{eq:Lipschitz_assumption}, \eqref{eq:bounded_assumption} and \eqref{eq:low-rank_assumption} hold \emph{globally}, because this greatly simplifies the analysis and the results. As in~\eqref{eq:assumption}, it is common to impose such assumptions only in a neighbourhood of the exact solution $A(t)$. When using randomized techniques, there is always a tiny but non-zero probability that the approximation leaves \emph{any} neighborhood. In Remark~\ref{remark: localized assumption}, we discuss how our analsyis can be modified to account for this effect, resulting in (slightly) weaker results.

\subsection{Low-rank approximability} \label{sec:lorapprox}

In this section, we discuss the relation between the assumption~\eqref{eq:assumption}
on the tangent space projection and the low-rank assumption~\eqref{eq:low-rank_assumption}.

First of all,~\eqref{eq:assumption} implies~\eqref{eq:low-rank_assumption}. To see this, suppose that $F$ satisfies~\eqref{eq:assumption}.
Then, by \cite[Lemma 1]{kieri2019projection}, there exists $h_0> 0$ such that \begin{equation*}
   \|\Phi^h_F(Y)-\Phi^h_{P_rF}(Y)\|_F\leq  \tilde{\epsilon}  \int^h_0e^{Ls}\,\mathrm{d}s\leq  e^{Lh}h\tilde{\epsilon}, \quad \forall\, 0\le h \le h_0.
\end{equation*} 
Because of $\Phi^h_{P_rF}(Y)\in \mathcal{M}_r$, this implies \begin{equation*}
     \|\Phi^h_F(Y)-[\![\Phi^h_F(Y)]\!]_r\|_F \leq \|\Phi^h_F(Y)-\Phi^h_{P_rF}(Y)\|_F\leq  e^{Lh}h\tilde{\epsilon}.
\end{equation*}
Hence,~\eqref{eq:low-rank_assumption} is satisfied with $C_M=e^{Lh_0}$ and $\epsilon=\tilde{\epsilon}$.

Assumption~\eqref{eq:low-rank_assumption} does not necessarily imply~\eqref{eq:assumption}, that is, the existence of a good low-rank approximation does not require the tangent space projection error to be small. In other words,
Assumption~\eqref{eq:low-rank_assumption} is weaker. This is demonstrated by the following example.
\begin{example} \label{example1}
Consider the rank-$2$ approximation of the differential equation
\begin{equation*}
    \dot{Y}(t)=F(Y(t))=\begin{pmatrix}
0 & 0 & 0\\
0 & 10 & 0\\
0& 0& -10
\end{pmatrix}Y(t)+\begin{pmatrix}
0 & 0 & 0\\
0 & 10^{-5}e^{-1} & 0\\
0& 0& 0
\end{pmatrix}, \qquad Y(0)=\begin{pmatrix}
1 & 0 & 0\\
0 & 0 & 0\\
0& 0& 10^{-6}e\\
\end{pmatrix}.
\end{equation*} We have that
\begin{equation*}
   \Phi^h_F(Y)= \begin{pmatrix}
1 & 0 & 0\\
0 & 10^{-6}(e^{10h-1} -e^{-1})& 0\\
0& 0& 10^{-6}e^{1-10h}\\
\end{pmatrix}
\end{equation*}
admits an excellent rank-$2$ approximation at $h = 1$:
$\|\Phi^1_F(Y)-[\![\Phi^1_F(Y)]\!]_2\|_F\leq 1.24\times 10^{-10}$. On the other hand,
the tangent space projection
\begin{equation*}
    \Phi^h_{PF}(Y)=\begin{pmatrix}
1 & 0 & 0\\
0 & 0& 0\\
0& 0& 10^{-6}e^{1-10h}\\ 
\end{pmatrix}.
\end{equation*}
results in a rank-$2$ matrix with a much larger error:
$\|\Phi^1_F(Y)-\Phi^1_{PF}(Y)\|_F\geq 0.008$.
\end{example}

 \subsection{Randomized low-rank approximation}
 \label{section: Randomized low-rank approximation}
 
 Conceptually, a rank-$r$ approximation is a map $\mathcal{R}:\mathbb{R}^{m\times n}\rightarrow \mathcal{M}_r$. When randomization is used, the map $\mathcal{R}$ is random, usually due to the use of random matrices for sketching. We measure the quality of $\mathcal{R}$ through moments, which will later be used to derive concentration inequalities. We say that $\mathcal{R}$ satisfies the moment assumption for $q\geq 1$  if 
\begin{equation}
    \label{eq: moment_assumption}\|\mathcal{R}(Z)-Z\|_{L_q}\leq C_\mathcal{R}\|Z-[\![Z]\!]_r\|_F,
\end{equation}
holds for fixed but arbitrary $Z\in \mathbb{R}^{m\times n}$, with a constant $C_\mathcal{R}$ being independent of $Z$.

Randomized low-rank approximations that utilize Gaussian random matrices for sketching usually satisfy the moment assumption~\eqref{eq: moment_assumption}. In the following, we will establish this fact for the generalized Nystr\"{o}m method~\cite{nakatsukasa2020fast,Tropp2017Practical}, which proceeds as follows. Given oversampling parameters $p,\ell \in \mathbb N$ and random matrices $\Omega\in \mathbb{R}^{n\times (r+p)}$, $\Psi\in \mathbb{R}^{m \times (r+p+\ell)}$,  generalized Nystr\"{o}m constructs a rank-$r$ approximation of $Z$ by first performing an oblique projection onto $\text{span}(Z\Omega)$ and then truncating to rank $r$:
\begin{equation}
\label{eq: Nystorm 1}
    Z\approx [\![ Z\Omega (\Psi^T Z\Omega)^\dagger\Psi^TZ]\!]_r:=\mathcal{N}(Z),
\end{equation}
where $(\cdot)^\dagger$ denotes the pseudoinverse.
If $\Psi^T Z \Omega\in \mathbb{R}^{(r+p+\ell)\times (r+p)}$ has full column rank, we have the equivalent expression 
\begin{equation}
\label{eq: Nystorm-2}
    \mathcal{N}(Z)= Q[\![(\Psi^T Q)^\dagger\Psi^TZ]\!]_r,
\end{equation}
where $Q$ is an orthonormal basis of span$(Z\Omega)$, computed by, e.g., a QR factorization, $Z\Omega = QR$.
For dense and unstructured matrices $\Omega$ and $\Psi$, 
computing the sketches $Z\Omega$ and $\Psi^TZ$ usually dominates the overall computational cost. In particular, this is true when 
$\Omega$ and $\Psi$ are Gaussian random matrices, i.e., their entries are independent standard normal Gaussian random variables.  The following theorem shows that  generalized Nystr\"{o}m satisfies the moment assumption \eqref{eq: moment_assumption} in this case.
\begin{theorem}
\label{theorem: ny}
    Suppose that $\Omega\in \mathbb{R}^{n\times (r+p)}$ and $\Psi\in \mathbb{R}^{m\times (r+p+\ell)}$ are independent standard Gaussian matrices with $p,\ell\geq 4$. Setting $q=\min\{p,\ell\}$, it holds for $Z\in \mathbb{R}^{m\times n}$ that
    \begin{equation*}
       \|\mathcal{N}(Z)-Z\|_{L_q}=\big(\mathbb{E}[\|\mathcal{N}(Z)-Z\|^q_F]\big)^{\frac{1}{q}}\leq C_{\mathcal{N}}\|Z-[\![ Z]\!]_r\|_F
    \end{equation*} with $C_{\mathcal{N}}=1+2\sqrt{(1+r+p)(1+r)}$.
\end{theorem}
\begin{proof}
By the triangle inequality
\begin{align*}
    \|\mathcal{N}(Z)-Z\|_{L_q} &\leq \|[\![Z\Omega (\Psi^T Z\Omega)^\dagger\Psi^TZ]\!]_r-Z\Omega (\Psi^T Z\Omega)^\dagger\Psi^TZ\|_{L_q}+\|Z\Omega (\Psi^T Z\Omega)^\dagger\Psi^TZ-Z\|_{L_q}
    \\&\leq \|[\![Z]\!]_r-Z\Omega (\Psi^T Z\Omega)^\dagger\Psi^TZ\|_{L_q}+\|Z\Omega (\Psi^T Z\Omega)^\dagger\Psi^TZ-Z\|_{L_q}
    \\&\leq \|[\![Z]\!]_r-Z\|_{L_q}+2\|Z\Omega (\Psi^T Z\Omega)^\dagger\Psi^TZ-Z\|_{L_q}.
\end{align*} To bound the second term,
we follow the proof of \cite[Theorem 11]{kressner2023randomized}. Considering an orthonormal basis $Q$ of $Z\Omega$, one obtains that
\begin{align*}
 \mathbb{E}[\|Z-Z\Omega (\Psi^T Z\Omega)^\dagger\Psi^TZ\|_F^{q}]&\leq (1+r+p)^{q/2}\mathbb{E}[\|(I-QQ^T)Z\|_F^{q}]\\&\leq (1+r+p)^{q/2}\left(\mathbb{E}^{p}[\|(I-QQ^T)Z\|_F]\right)^q\\&=(1+r+p)^{q/2}\left(\mathbb{E}^{p/2}[\|(I-QQ^T)Z\|_F^{2}]\right)^{q/2}\\ &\leq (1+r+p)^{q/2}\left((1+r)\|Z-[\![ Z]\!]_r\|_F^2\right)^{q/2}.
\end{align*} Therefore,
$\|\mathcal{N}(Z)-Z\|_{L_q} \leq \|Z-[\![ Z]\!]_r\|_F+  2\sqrt{(1+r+p)(1+r)}\|Z-[\![ Z]\!]_r\|_F = C_{\mathcal{N}}\|Z-[\![ Z]\!]_r\|_F$.
\end{proof} 

To keep our developments concrete, we will always use generalized Nystr\"{o}m 
instead of an abstract randomized method $\mathcal{R}$ in the rest of the paper. However, the theoretical results remain valid for any $\mathcal{R}$  satisfying~\eqref{eq: moment_assumption}. 

\subsection{Idealized randomized projection method}

Following~\cite[Ch. IV.4]{Hairer2010Geometric}, a rank-$r$ approximation to the solution $A((i+1)h)$ is obtained by combining exact integration with rank-$r$ truncation
\begin{equation*}
    Y_{i+1}=[\![ \Phi^h_F(Y_i) ]\!]_r, \qquad Y_0 = [\![ A_0 ]\!]_r.
\end{equation*} 
This is an idealized integrator because the exact flow $\Phi^h_F$ still needs to be approximated in order to obtain a practical method. Replacing rank-$r$ truncation by the 
generalized Nystr\"{o}m method one gets 
the \emph{idealized randomized projection method}
\begin{equation}
\label{eq: iteration}
    Y_{i+1}=\mathcal{N}_i(\Phi^h_F(Y_i))=[\![ \Phi^h_F(Y_i)\Omega_i(\Psi_i\Phi^h_F(Y_i)\Omega_i)^\dagger\Psi_i^T\Phi^h_F(Y_i)]\!]_r. \qquad Y_0 = \mathcal{N}_0(A_0).
\end{equation}
The subscript $i$ of $\mathcal{N}$ is used to emphasize that the generalized Nystr\"{o}m method is used with different (independent) $\Omega_i\in \mathbb{R}^{n\times (r+p)}$, $\Psi_i\in \mathbb{R}^{m\times (r+p+\ell)}$ in every time step.

\subsubsection{Error analysis}

In this section, we provide an error analysis of the idealized method~\eqref{eq: iteration} when $\Omega_i$, $\Psi_i$ are independent standard Gaussian matrices. The following theorem provides a bound on the $L_q$ norm of the error.
The proof basically follows from the proof of \cite[Theorem 2]{kieri2019projection}, with the Frobenius norm replaced by the $L_q$ norm. It is included for convenience, because similar arguments will be used again below. 
\begin{theorem}
\label{Thm: idealised lq}
     With the assumptions stated in Section~\ref{sec:lorapprox} and assuming {$\|[\![ A_0 ]\!]_r-A_0\|_F\leq \delta$} holds for the initial data, the method \eqref{eq: iteration} with independent standard Gaussian matrices $\Omega_i\in \mathbb{R}^{n\times (r+p)}$, $\Psi_i\in \mathbb{R}^{m\times (r+p+\ell)}$ and oversampling parameters $p,\ell\geq 4$ satisfies the error estimate
    \begin{equation*}
       \|Y_N-A(Nh)\|_{L_q}\leq C(\delta+\epsilon)
    \end{equation*} for $q=\min\{p,\ell\}$  on a finite time-interval $0\leq Nh\leq T$ for every $0<h\leq h_0$. The constant $C$ only depends on $L,T,h_0$ $C_M$, and $C_{\mathcal{N}}$.
\end{theorem}
   \begin{proof}
   We first note that $\mathcal{N}_i$ is stochastically independent of $\Phi^h_F(Y_i)$. By the the law of
total expectation, Theorem \ref{theorem: ny} and Assumption~\eqref{eq:low-rank_assumption}, we get \begin{align}
    \label{eq: local error}
        \|\mathcal{N}_i(\Phi^h_F(Y_i))-\Phi^h_F(Y_i)\|_{L_q}& =\Big(\mathbb{E}\Big[\mathbb{E}[\|\mathcal{N}_i(\Phi^h_F(Y_i))-\Phi^h_F(Y_i)\|^{q}_F |Y_i]\Big]\Big)^{1/q}\nonumber\\&\leq\Big( \mathbb{E}[C_{\mathcal{N}}^q\|\Phi^h_F(Y_i)-[\![\Phi^h_F(Y_i)]\!]_r\|^q_F]\Big)^{1/q} \leq C_{\mathcal{N}}C_M\epsilon h.
    \end{align}  To bound the $L_q$ norm of the global error, we follow the proof of \cite[Theorem 2]{kieri2019projection} and use a telescoping sum:
    \begin{equation*}
        \|Y_N-A(Nh)\|_{L_q}=\Big\|\sum^N_{i=1}(\Phi^{(N-i)h}_{F}(Y_i)-\Phi^{(N-i+1)h}_{F}(Y_{i-1}))+\Phi^{Nh}_{F}(Y_0)-\Phi^{Nh}_{F}(A_{0})\Big\|_{L_q}\leq \sum^N_{i=0}E_i,
    \end{equation*} where we define 
    \begin{align*}
        E_i&=\|\Phi^{(N-i)h}_{F}(Y_i)-\Phi^{(N-i)h}_{F}(\Phi^{h}_{F}(Y_{i-1}))\|_{L_q},\quad i=1,\cdots,N, \\E_0&=\|\Phi^{Nh}_{F}(Y_0)-\Phi^{Nh}_{F}(A_{0})\|_{L_q}.
    \end{align*} The Lipschitz continuity of $F$ implies
    $\|\Phi^{t}_{F}(X)-\Phi^{t}_{F}(Y)\|_F\leq e^{Lt}\|X-Y\|_F$.
    Therefore, by the law of
total expectation and \eqref{eq: local error},
    \begin{align}
    \label{eq: bound of local error}
        E_i&=\Big(\mathbb{E}\Big[\mathbb{E}[\|\Phi^{(N-i)h}_{F}(Y_i)-\Phi^{(N-i)h}_{F}(\Phi^{h}_{F}(Y_{i-1}))\|^{q}_F| Y_{i-1}]\Big]\Big)^\frac{1}{q}\nonumber\\&\leq e^{Lh(N-i)}\Big(\mathbb{E}\Big[\mathbb{E}[\|Y_i-\Phi^{h}_{F}(Y_{i-1})\|^{q}_F| Y_{i-1}]\Big]\Big)^\frac{1}{q}\nonumber\\&= e^{Lh(N-i)}\Big(\mathbb{E}\Big[\mathbb{E}[\|\mathcal{N}_{i-1}(\Phi^h_F(Y_{i-1}))-\Phi^{h}_{F}(Y_{i-1})\|^{q}_F| Y_{i-1}]\Big]\Big)^\frac{1}{q}\nonumber\\&\leq C_{\mathcal{N}}\cdot C_M e^{Lh(N-i)}\epsilon h.
    \end{align} In summary,
    \begin{equation*}
         \|Y_N-A(Nh)\|_{L_q}\leq C_{\mathcal{N}}e^{LNh}\delta+C_{\mathcal{N}}\cdot C_M\epsilon\sum he^{Lh(N-i)},
    \end{equation*} which yields the desired result by bounding the sum by an integral, as in \cite[P.80]{kieri2019projection}.
\end{proof}
The Markov inequality turns the moment bound of Theorem \ref{Thm: idealised lq} into tail bounds for the approximation error.
\begin{corollary}
\label{corollary: idealised_method_prob}
  With the assumptions and notation stated in Theorem \ref{Thm: idealised lq},
  the error estimate
  \begin{equation*}
        \|Y_N-A(Nh)\|_F\leq C\eta(\delta+\epsilon),
\end{equation*}
holds for any $\eta \geq 1$
with probability at least $1-\eta^{-q}$.
%
\end{corollary}
\begin{proof}
By Markov's inequality and Theorem \ref{Thm: idealised lq},
\begin{equation*}
    \Pr\{ \|Y_N-A(Nh)\|_F\geq C\eta(\delta+\epsilon)\}\leq \left(\frac{[\mathbb{E}\|Y_N-A(Nh)\|_F^q]^{1/q}}{C\eta(\delta+\epsilon)}\right)^q\leq\frac{1}{\eta^q}.
\end{equation*}
\end{proof}

Corollary~\ref{corollary: idealised_method_prob} states that the generalized Nystr\"{o}m method produces an error on the level of $\epsilon+\delta$ with high probability.
Under the assumptions stated in Section~\ref{sec:lorapprox},
the same type of error bound is obtained when using exact rank-$r$ truncations.


\section{Randomized low-rank Runge-Kutta methods}
\label{Section: 3}

To turn the idealized projection method~\eqref{eq: iteration} into a practical method, we need to combine it with a time-integration method, e.g., a RK method.
However, directly replacing the exact flow $\Phi_F$ by a RK method will result in high ranks in the intermediate stages. To mitigate this issue, we also apply the generalized Nystr\"{o}m method to these intermediate stages. To make this idea specific, let us consider a general explicit Runge-Kutta method with $s$ stages applied to the matrix differential equation~\eqref{ode}: \begin{equation}
   \label{eq: standard RK}
    \begin{aligned}
        \tilde{Z}_j&=A_i+h\sum^{j-1}_{l=1}a_{jl}F(\tilde{Z}_l),\quad j=1,\ldots,s,\\ 
        A_{i+1}&=A_i+h\sum^{s}_{j=1}b_jF(\tilde{Z}_j).  \end{aligned} 
        \end{equation}
Performing generalized Nystr\"{o}m in the intermediate stages yields our \emph{Randomized low-rank RK method}:
\begin{equation}
\label{rand_rk}
    \begin{aligned}
        Z_j&=Y_i+h\sum^{j-1}_{l=1}a_{jl}F({\mathcal{N}_l}(Z_l)),\quad j=1,\ldots,s,
        \\
        Y_{i+1}&=\mathcal N_{s+1}\Big(Y_i+h\sum^{s}_{j=1}b_jF({\mathcal{N}_j}(Z_j))\Big). \end{aligned} 
\end{equation}
Note that the index of $\mathcal N$ is now used to emphasize the use of different (independent)
random matrices $\Omega_j, \Psi_j$ for different stages. Across different time steps, the random matrices are, of course, also independently drawn.
In~\eqref{rand_rk}, the rank of $Z_j$ increases linearly with respect to $j$. However, there is no need to construct and store $Z_j$ explicitly. For all subsequent purposes, we only need access to the Nystr\"om approximation of $Z_j$, and thus it suffices to compute and store the sketches $Z_j \Omega_j$ and $\Psi_j^TZ_j$. In fact, the method~\eqref{rand_rk} is equivalent to 
\begin{equation}
\label{rand_rk_2}
    \begin{aligned}
    &\begin{cases}
        \displaystyle Z_j\Omega_j=Y_i\Omega_j+h\sum^{j-1}_{l=1}a_{jl}F\left({\mathcal{N}_l}(Z_l)\right)\Omega_j
        \\
        \displaystyle \Psi^T_jZ_j=\Psi_j^TY_i+h\sum^{j-1}_{l=1}a_{jl}\Psi^T_jF\left({\mathcal{N}_l}(Z_l)\right),
    \end{cases}\quad j=1,\ldots,s,
    \\&\quad Y_{i+1}=\mathcal{N}_{s+1}\Big(Y_i+h\sum^{s}_{j=1}b_jF\left({\mathcal{N}_j}(Z_j)\right)\Big).   \end{aligned} 
\end{equation}
We will use the expression~\eqref{eq: Nystorm-2} to evaluate $\mathcal{N}_j(Z_j)$ and, for this purposes, only needs $\Psi_j^TZ_j$, $\Psi_j$, and $Z_j\Omega_j$ (or, rather, the orthogonal factor of its QR decomposition).

Algorithm \ref{alg:rrk-s} contains the pseudo-code of the Randomized low-rank RK method. It closely follows~\eqref{rand_rk_2}, except that we also precompute the sketches of $F([\![ \hat{Z_j}(\Psi^T_j \hat{Z_j})^\dagger\tilde{Z}_j]\!]_r)$ as they are needed in subsequent stages.
\begin{algorithm}[h]
\caption{Randomized low-rank Runge-Kutta method with $s$ stages}\label{alg:rrk-s}
\begin{algorithmic}
\Require Differential equation~\eqref{ode} defined by $F$ and initial condition $A_0\in \mathbb{R}^{m \times n}$. Target rank $r$, oversampling parameters $p,\ell$, step size $h>0$, number of steps $N\geq0$.
\Ensure Approximation $Y_N\in \mathcal{M}_r$ of $A(Nh)$.
\State Draw ind. random matrices $\Omega\in \mathbb{R}^{n \times (r+p)}$, $\Psi\in \mathbb{R}^{m \times (r+p+\ell)}$.
\State $\hat{Y}_0=A_0\Omega$, $\tilde{Y}_0=\Psi^TA_0$
\For{$i=0,\ldots, N-1$}
\State Draw ind. random matrices $\Omega_j\in \mathbb{R}^{n \times (r+p)}$, $\Psi_j\in \mathbb{R}^{m \times (r+p+\ell)}$ for $j=1,2,\ldots, (s+1)$.
\For{$j=1,\ldots, s$}
\State $\hat{Z}_j= [\![ \hat{Y}_{i}(\Psi^T\hat{Y}_{i})^\dagger\tilde{Y}_{i}]\!]_r\Omega_j+h\sum^{j-1}_{l=1}a_{jl}\hat{K}_{lj}$ \Comment{Evaluate $ [\![ \hat{Y}_{i}(\Psi^T\hat{Y}_{i})^\dagger\tilde{Y}_{i}]\!]_r$ using \eqref{eq: Nystorm-2}}
\State $\tilde{Z}_j= \Psi^T_j [\![ \hat{Y}_{i}(\Psi^T\hat{Y}_{i})^\dagger\tilde{Y}_{i}]\!]_r+h\sum^{j-1}_{l=1}a_{jl}\tilde{K}_{lj}$
\State $F_j=F([\![ \hat{Z_j}(\Psi^T_j \hat{Z_j})^\dagger\tilde{Z}_j]\!]_r)$\Comment{Evaluate $[\![ \hat{Z_j}(\Psi^T_j \hat{Z_j})^\dagger\tilde{Z}_j]\!]_r$ using~\eqref{eq: Nystorm-2}}
\For{$q=j+1,\ldots,s+1$} \Comment{Pre-compute sketches of $F$ for other stages}
\State $\hat{K}_{jq}=F_j\Omega_q$ 
\State $\tilde{K}_{jq}=\Psi^T_qF_j$ 
\EndFor

\EndFor
\State $\hat{Y}_{i+1}=Y_i\Omega_{s+1}+h\sum^{s}_{j=1}b_j\hat{K}_{j,s+1}$
\State $\tilde{Y}_{i+1}=\Psi^T_{s+1}Y_i+h\sum^{s}_{j=1}b_j\tilde{K}_{j,s+1}$
\State Set $\Psi=\Psi_{s+1}$
\EndFor
\State Return $Y_{N}=[\![ \hat{Y}_{N}(\Psi\hat{Y}_{N})^\dagger\tilde{Y}_{N}]\!]_r.$ \Comment{Compute $ [\![ \hat{Y}_{N}(\Psi^T\hat{Y}_{N})^\dagger\tilde{Y}_{N}]\!]_r$ using~\eqref{eq: Nystorm-2}}
\end{algorithmic}
\end{algorithm}

\subsection{Implementation aspects and cost} \label{sec:implementation}

We now consider the efficient implementation and cost of a time step performed by Algorithm~\ref{alg:rrk-s}. To simplify the discussion, we assume $m=n$ and $p=\ell=\mathcal{O}(r) \ll n$. We let $c_n$ denote the cost of multiplying a vector of length $n$ with $\Omega$ or $\Psi$. In the worst case, when $\Omega, \Psi$ are unstructured dense random matrices, $c_n=\mathcal{O}(nr)$. The use of structured random matrices can lead to lower $c_n$. For example, using the Subsampled Randomized Fourier Transform (SRFT)~\cite{Halko2011Finding} for sketching reduces $c_n$ to $\mathcal{O}(n\log(r))$. 

Every (large) $n\times n$ matrix occurring in the algorithm is represented in factored form $U\Sigma V^T\in\mathbb{R}^{n\times n}$,
where $U, V$ are tall matrices (not necessarily orthonormal) and $\Sigma$ is a small square matrix (not necessarily diagonal) of size equal to the rank of the matrix.

The evaluation of $[\![ \hat{Z_j}(\Psi^T_j \hat{Z_j})^\dagger\tilde{Z}_j]\!]_r$ in Algorithm~\ref{alg:rrk-s} requires $\mathcal{O}(nr^2+jnr)$ operations for computing $\hat{Z_j}$, $\tilde{Z_j}$ and another $\mathcal{O}(nr^2)$ operations for computing the factored form of the matrix.

The cost of applying $F$ to a rank-$r$ matrix and obtaining a factored form of the result strongly depends on the nature of $F$. We will denote this cost by $c_F$ and the resulting rank by $r_F$. In some cases (see Section~\ref{sec:lyapunovexample} for an example), $c_F = \mathcal O( nr )$ and $r_F = \mathcal O(r)$. In cases that lead to large $r_F$, the use of random sketches gives flexibility to exploit structure. For example, consider the case that $F(A)$ contains Hadamard products of the matrix $A$ with itself, originating from, e.g., quadratic nonlinearities in the underlying partial differential equation. Then although the rank of the matrix $F_j$ in Algorithm~\ref{alg:rrk-s} is much larger than $r$, its factored representation has rich Kronecker product structure, which can be exploited when sketching $F_j$ with Khatri-Rao products of random matrices~\cite{bujanovic2024subspace,KressnerPerisa17}. As another example, when $F$ has block structure, one can use the Block SRFT \cite{Balabanov2023Block} to benefit from parallel computing. In Section \ref{section 4.1.1} , we test numerically the possibility to speed up the computation of $\hat{K}_{jq}$ and $\tilde{K}_{jq}$ by using the same random matrices. Although there is little justification, this variant appears to lead to an accuracy comparable to Algorithm~\ref{alg:rrk-s}. When not using any of the tricks mentioned above, the 
computation of each $\hat{K}_{jq}$ and $\tilde{K}_{jq}$ requires $\mathcal{O}(r_Fc_n+r_F^2r+nr_Fr) = \mathcal{O}(r_Fc_n+nr_Fr)$ operations.

In summary, the complexity of the $j$th stage is
$\mathcal{O}(jnr+nr^2+c_F+(s+1-j)(r_Fc_n+nr_Fr)),$ and thus one step of the Randomized RK method has a total complexity of
\begin{equation} \label{eq:complexityrandrk} \mathcal{O}(s^2(r_Fc_n+nr_Fr)+snr^2+sc_F).\end{equation} It can be seen that the computation of $\hat{K}_{jq}$, $\tilde{K}_{jq}$ is a dominating part of the cost. 

The linearity of the sketches with respect to the data (sometimes also called a streaming property), makes the generalized Nystr\"{o}m method a preferred choice for the randomized low-rank approximation in Algorithm~\ref{alg:rrk-s}. It allows one to only store and work with sketches of $F_j$, $j=1,\ldots, s$ when computing subsequent stages. In contrast, the randomized SVD uses an orthogonal projection that is not linear in the data and, in turn, the full rank-$r_F$ factorizations of the potentially high-rank matrices $F_j$ need to be stored. Another disadvantages of the randomized SVD is that its direct application to the sums appearing in the Runge-Kutta stages is quite expensive, requiring $\mathcal{O}(ns^2r_F^2)$ operations.
    
\subsection{Comparison to Projected Runge-Kutta method} \label{sec:compprk}

{The Projected Runge-Kutta method (PRK) from~\cite{kieri2019projection} is closely related to our proposed method~\eqref{rand_rk}. It proceeds by performing the time stepping \begin{equation}
\label{eq: PRK}
    \begin{aligned}
        Z_j&=Y_i+h\sum^{j-1}_{l=1}a_{jl}{P_r}(\mathcal{R}(Z_l))F(\mathcal{R}(Z_l)),\quad j=1,\ldots,s,
        \\
        Y_{i+1}&=\mathcal{R}\Big(Y_i+h\sum^{s}_{j=1}b_jP_r(\mathcal{R}(Z_j))F(\mathcal{R}(Z_j))\Big).   \end{aligned} 
\end{equation}
Here, $\mathcal R$ denotes a retraction to the manifold $\mathcal M_r$ and a common choice is the truncated SVD.
The tangent space projection $P_r$ in~\eqref{eq: PRK} reduces the rank of $F(\mathcal{R}(Z_l))$ to $2r$, which can make the subsequent application of $\mathcal R$ signficantly cheaper.
Our method uses sketching instead of tangent space projection, which has two potential advantages: (1) As discussed in Section~\ref{sec:lorapprox}, the tangent space projection could introduce significant error, even when the solution admits a good rank-$r$ approximation. In this situation, we expect our method~\eqref{rand_rk} to be more accurate. This expectation is confirmed by the numerical experiments in Section~\ref{sect: Numerical Experiments}. (2) As discussed in Section~\ref{sec:implementation} above, sketching gives additional flexibility in the choice of random matrices to exploit structure in $F$ applied to a rank-$r$ matrix. Tangent space projection does not offer this flexibility.

One iteration of $\eqref{eq: PRK}$ needs to apply retraction to the matrices $Z_j$, which have rank at most $2jr$ for $j=1\ldots s$. When using the truncated SVD to perform retraction, the cost
is dominated by the QR factorizations of the $n\times 2jr$ factors of $Z_j$. In summary, the total complexity for one time step of PRK~\eqref{eq: PRK} is
\begin{equation*}
    \mathcal{O}(\underbrace{s^3nr^2}_{\text{total retraction cost}}+\underbrace{snr_Fr}_{\text{$s \times$ application of $P_r$}}+\underbrace{sc_F}_{\text{$s \times$ evaluation of $F$}})
\end{equation*}
Compared to~\eqref{eq:complexityrandrk}, we see that the term $s^2nr_Fr$ is reduced to $snr_Fr$, which only becomes relevant when $r_F> sr$. This is because PRK can reuse computations related to tangent space projections across different stages, while the randomized low-rank RK method uses different sketches for every stage and can thus not reuse computations. As mentioned above, this issue can be mitigated by reusing random matrices for sketching, at the expense of theoretical justification.
However, as $s$ is typically very small (e.g., $1$, $2$ or $4$), this issue may not be too relevant in practice.

\subsection{Error analysis}
\label{sect:3.3}

With high probability, our method achieves at least the same qualitative error behavior that has been established for PRK. To see this, we follow~\cite{kieri2019projection} and consider the stage orders $\gamma_1,\ldots \gamma_s$, which are defined as the local errors of the stages $\tilde{Z_j}$ in the standard RK method~\eqref{eq: standard RK}: For every $h\leq h_0$,
\begin{equation} \label{eq:stageorder}
    \|\tilde{Z}_j-\phi^{c_jh}_F(A_i)\|_F\leq C_Lh^{\gamma_j+1},\quad j=1, \ldots,s,
\end{equation} with $c_j = a_{j1} +a_{j2} + \cdots + a_{j,j-1}$.
We then obtain the following result, which corresponds to Theorem 6 in~\cite{kieri2019projection}.
\begin{theorem}
\label{Thm: stage order}
    Consider the randomized low-rank RK method~\eqref{rand_rk} utilizing independent standard Gaussian matrices, oversampling parameters $p,\ell\geq 4$, and an explicit $s$-stage RK method of order $\tau$ with stage orders $\gamma_1\leq \gamma_2\leq \cdots \leq \gamma_s$. Denote 
    \begin{equation*}
        \gamma=\begin{cases}
            \min(\tau,\gamma_2+1) & \text{if } b_2\neq 0,\\
            \min(\tau,\gamma_3+1,\gamma_2+2) & \text{if } b_2= 0. 
        \end{cases}
    \end{equation*}Then with the assumptions stated in Theorem~\ref{Thm: idealised lq}, the global error is bounded for $q=\min\{p,\ell\}$ by 
    \begin{equation*}
       \|Y_N-A(Nh)\|_{L_q}\leq C(\delta+\epsilon+h^\gamma)
    \end{equation*} on the finite time interval $0\leq Nh\leq T$, for all $h\leq h_0$. The constant $C$ depends only on $L$,$T$, $h_0$,$s$, $C_L$,$C_{\mathcal{N}}$, $\max_{ij}|a_{ij}|$ and $\max_{i}|b_{i}|$. In particular, for any $\eta \geq 1$, it holds for fixed $h$ and $N$ that
\begin{equation*}
    \Pr\big\{ \|Y_N-A(Nh)\|_F\geq C\eta(\epsilon+h^\gamma+\delta) \big\}\leq\frac{1}{\eta^q}.
\end{equation*}
\end{theorem}
 \begin{proof}
The result of this theorem essentially follows from replacing the Frobenius norm in the proof of~\cite[Theorem 6]{kieri2019projection} by the $L_q$ norm and performing some additional minor modifications. For completeness, we have included the proof in the appendix.
\end{proof}

Theorem~\ref{Thm: stage order} above shows that the randomized RK methods with the coefficients given by the following RK methods of order 1,2 and 3, enjoy the usual convergence order up to $\mathcal{O}(\epsilon)$:
\begin{itemize}
    \item RK 1 (Euler): $b_1=1,$
    \item RK 2 (Heun's method): $a_{21}=1$, $b_1=b_2=\frac{1}{2},$
    \item RK 3 (Heun's third-order method): $a_{21}=a_{32}=\frac{1}{2}$, $a_{43}=1,$ $b_1=b_4=\frac{1}{6}$, $b_2=b_3=\frac{1}{3}$.
\end{itemize}
Unfortunately, for RK 4 we only obtain order $2$ from Theorem~\ref{Thm: stage order}. Section~\ref{sec:rk4} investigates this combination further.

\begin{remark}
\label{remark: localized assumption}
As already noted, the global nature of the three assumptions in Section~\ref{section: assumptions} can be quite limiting.
If we modify these assumptions~\eqref{eq:Lipschitz_assumption},~\eqref{eq:bounded_assumption} and~\eqref{eq:low-rank_assumption} such that they only hold in a neighborhood of the exact solution $A(t)$ for $0 \leq t\leq T$, we need to additionally ensure that $Y_i$ and the intermediate stages remain in the neighborhood. Due to the presence of randomness, this complicates the analysis and yields slightly worse results. In the following, we sketch which modifications need to be performed in order to localize the assumptions.

Suppose we want to ensure $E_i:=\|Y_i-A(ih)\|_F\leq M$ for every $i=1,\ldots N-1$ to use the properties \eqref{eq:Lipschitz_assumption}, \eqref{eq:bounded_assumption} and \eqref{eq:low-rank_assumption} locally. (We simplify the discussion by ignoring the probability that the intermediate stages leave the neighborhood; it is easy to adapt the approach outlined here by additionally requiring $\|Z_j-\tilde{Z}_j\|_F\leq M'$ for $j=1,\ldots, s$, which will yield the same convergence rate with a different constant and a somewhat higher failure probability, increased by a factor $s$.) To proceed, we can utilize the failure probability estimate of Theorem~\ref{Thm: stage order} to conclude
\begin{equation*}
    \Pr\big\{E_i> M | \cap^{i-1}_{j=0} \{E_{j}\leq M\}\big\}\leq\left(\frac{C(\epsilon+h^\gamma+\delta)}{M}\right)^{q}. 
\end{equation*} Using conditional probability and Bernoulli's inequality,
\begin{align*}
    \Pr\big\{ \cap^{N-1}_{j=0} \{E_{j}\leq M\}\big\}&=\Pr\big\{E_{N-1}\leq M |\cap^{N-2}_{j=0} \{E_{j}\leq M\}\big\}\Pr\big\{ \cap^{N-2}_{j=0} \{E_{j}\leq M\}\big\}\\&\geq \Big(1-\left(\frac{C(\epsilon+h^\gamma+\delta)}{M}\right)^{q}\Big)\Pr\{ \cap^{N-2}_{j=0} \{E_{j}\leq M\}\}
    \\&\geq \Big(1-\left(\frac{C(\epsilon+h^\gamma+\delta)}{M}\right)^{q}\Big)^{(N-1)}
    \\&\geq 1-(N-1)\left(\frac{C(\epsilon+h^\gamma+\delta)}{M}\right)^{q}.
\end{align*}
Similarly, we have \begin{equation*}
    \Pr\{ \|Y_N-A(Nh)\|_F\geq C\eta(\epsilon+h^\gamma+\delta)\}\leq \frac{1}{\eta^q}+(N-1)\left(\frac{C(\epsilon+h^\gamma+\delta)}{M}\right)^{q}. 
\end{equation*}  The additional second term is not large when $M$ is large and/or $\epsilon+h^\gamma+\delta$ is small. To further quantify how this term affects the order of convergence, we substitute $\eta=\frac{M}{h^{1/q}}$ and assume $C(\epsilon+h^\gamma+\delta)\leq 1$, $h\leq 1$, leading to
\begin{equation}
\label{eq: total prob}
    \Pr\left\{ \|Y_N-A(Nh)\|_F\geq {CM}(\epsilon+h^\gamma+\delta)h^{-\frac{1}{q}}\right\}\leq \frac{h}{M^q}+(N-1)\left(\frac{C(\epsilon+h^\gamma+\delta)}{M}\right)^{q}\leq \frac{N}{M^q}.
\end{equation}
The presence of the additional factor $h^{-\frac{1}{q}}$ reduces the convergence order by ${\frac{1}{q}}$. Because of $q = \min\{p,\ell\}$, even modest choices of the oversampling parameters $p,\ell$ imply that this potential order loss is negligible.
´\end{remark}

\subsection{Error analysis of randomized low-rank Runge-Kutta 4 method}\label{sec:rk4}

Although our numerical experiments indicate convergence order $4$ for the randomized RK method based on RK 4, it appears to be difficult to establish this order theoretically.
In the following, we establish order $4$ when the intermediate stages are oversampled. Let us emphasize that this oversampling is only performed for theoretical purposes;
in practice, it does not seem to be needed.

Concretely, we plug the coefficients of the classical RK 4 method into~\eqref{rand_rk} and oversample the intermediate stage as follows:
\begin{align}
\label{eq: modified scheme}
    \hat{Z}_1&=\mathcal{N}_1^{15r}({Y}_i)\nonumber\\
    \hat{Z}_2&=\mathcal{N}_2^{15r}({Y}_i+\frac{h}{2}F(\hat{Z}_1))\nonumber\\
    \hat{Z}_3&=\mathcal{N}_3^{15r}({Y}_i+\frac{h}{2}F(\hat{Z}_2))\nonumber\\
    \hat{Z}_4&=\mathcal{N}_4^{15r}({Y}_i+{h}F(\hat{Z}_3))\nonumber\\
  {Y}_{i+1}&=\mathcal{N}_5\left({Y}_i+\frac{h}{6}(F(\hat{Z}_1)+2F(\hat{Z}_2)+2F(\hat{Z}_3)+F(\hat{Z}_4))\right).
\end{align} Here, $\mathcal{N}_5$ refers to the usual generalized Nystr\"{o}m method~\eqref{eq: Nystorm 1} with target rank $r$, while
$\mathcal{N}_i^{15r}$, for $i =1, \ldots,4$, refers to the generalized Nystr\"{o}m method with the increased target rank $15r$. 
As the method leaves $\mathcal{M}_r$, we need to impose a stronger assumption on low-rank approximability. For $h\leq h_0$, we assume that
 \begin{equation}
  \label{eq:low-rank_assumption_adaptive}
    \|\Phi^h_F(Y)-[\![\Phi^h_F(Y)]\!]_{k}\|_{F}\leq C_Mh\epsilon, \quad \forall\,Y\in \mathcal{M}_{k}, \quad \forall\,r\leq k\leq15r.
 \end{equation}

We start our analysis of~\eqref{eq: modified scheme} with a result on the low-rank approximability of $F$ implied by~\eqref{eq:low-rank_assumption_adaptive}. 
\begin{lemma}
\label{lemma: low_rank_F}
    Assuming that $F$ satisfies $\eqref{eq:low-rank_assumption_adaptive}$, let $k$ be any integer such that $r\leq k\leq 15r$. Then 
    \begin{equation*}
        \|F(Y)-[\![F(Y)]\!]_{2k}\|_{F}\leq C_M\epsilon, \quad \forall\, Y\in \mathcal{M}_{k}.
    \end{equation*}
\end{lemma}\begin{proof}
By the definition of the flow, $F(Y)$ is the time derivative of $\Phi^t_F(Y)$. Thus, for any $\gamma > 0$, there exists $h>0$ such that
%
%
%
    \begin{equation*}
        \Big\|\frac{\Phi^{h}_F(Y)-\Phi^0_F(Y)}{h}-F(Y)\Big\|_{F}\leq \gamma.
    \end{equation*}Because $[\![\Phi^{h_j}_F(Y)]\!]_{k}-\Phi^0_F(Y)$ has rank at most $2k$, it follows that 
    \begin{align*}
         \|F(Y)-[\![F(Y)]\!]_{2k}\|_{F}&\leq \Big\|F(Y)-\frac{[\![\Phi^{h_j}_F(Y)]\!]_{k}-\Phi^0_F(Y)}{h_j}\Big\|_{F}\\
         &\leq  \Big\|F(Y)-\frac{\Phi^{h_j}_F(Y)-\Phi^0_F(Y)}{h_j}\Big\|_{F}+\Big\|\frac{\Phi^{h_j}_F(Y)-[\![\Phi^{h_j}_F(Y)]\!]_{k}}{h_j}\Big\|_{F}
         \\&\leq \gamma+C_M \epsilon.
    \end{align*}
    The result of the lemma is obtained by taking $\gamma \rightarrow 0$.
\end{proof}
The following auxiliary result helps to bound the local error of~\eqref{eq: modified scheme}.
\begin{lemma}
\label{lemma: low-rank bound}
Let $Z\in \mathcal{M}_r, X\in \mathbb{R}^{m\times n}$, $\alpha\geq 0$, and $k\leq 14r$. Then
    \begin{equation*}
        \|Z+\alpha X-[\![Z+\alpha X]\!]_{15r}\|_F \leq \alpha \|X-[\![B]\!]_{k}\|_F
    \end{equation*}
    holds for any $B\in \mathbb{R}^{m\times n}$.
\end{lemma}
\begin{proof} Using that $Z+\alpha [\![B]\!]_{k}$ has rank at most $15r$, the result follows from
\[
            \|Z+\alpha X-[\![ Z+\alpha X]\!]_{15r}\|_F \leq \|Z+\alpha X-[\![Z+\alpha[\![B]\!]_{k}]\!]_{15r}\|_F
             =\alpha\|X-[\![B]\!]_{k}\|_F.
 \]
 \end{proof}
We are now in the position to establish a local error estimate for~\eqref{eq: modified scheme}.
\begin{lemma}
\label{lemma: local error of Yhat}
Suppose that the assumptions stated in Section~\ref{section: assumptions}
and~\eqref{eq:low-rank_assumption_adaptive} hold. Given $Y_i\in \mathcal{M}_r$, one step of the method~\eqref{eq: modified scheme} with independent standard Gaussian matrices and oversampling parameters $p,\ell\geq 4$ satisfies the local error estimate
    \begin{equation*}
        \|\Phi^{h}_F({Y}_{i})-{Y}_{i+1}\|_{L_q}\leq C(h\epsilon+h^5),
    \end{equation*}for $q=\min\{p,\ell\}$ and all $0<h\leq h_0$. The constant $C$ depends only on $L,T,h_0$ $C_M$, and $C_{\mathcal{N}}$.
\end{lemma}
\begin{proof}
The proof proceeds by bounding the moments of the differences between the stages $\hat{Z}_j$ of~\eqref{eq: modified scheme} and the stages $\tilde{Z}_j$ of the classic RK4 applied to ${Y}_i$, as defined in \eqref{eq: standard RK}.
For $j = 1$, ${Y}_i\in \mathcal{M}_r$ implies
    \begin{equation*}
        \|\hat{Z_1}-\tilde{Z_1}\|_{L_q}=\|\mathcal{N}_1^{15r}({Y}_i)-{Y}_i\|_{L_q}=0.
    \end{equation*}
For $j = 2$, we set $Y_{i,1}:={Y}_i+\frac{h}{2}F(\hat{Z}_1)$
    \begin{align}
         \|\hat{Z_2}-\tilde{Z_2}\|_{L_q}&=\Big\|\mathcal{N}_2^{15r}(Y_{i,1})-{Y}_i-\frac{h}{2}F(\tilde{Z}_1)\Big\|_{L_q}\nonumber
        \\&\leq \Big\|\mathcal{N}_2^{15r}(Y_{i,1})-\mathcal{N}_2^{15r}({Y}_i+\frac{h}{2}[\![F(\tilde{Z}_1)]\!]_{2r})\Big\|_{L_q} \nonumber 
        \\&\quad +\Big\|\mathcal{N}_2^{15r}({Y}_i+\frac{h}{2}[\![F(\tilde{Z}_1)]\!]_{2r})-{Y}_i-\frac{h}{2}F(\tilde{Z}_1)\Big\|_{L_q}\nonumber
        \\&=\Big\|\mathcal{N}_2^{15r}(Y_{i,1})-\mathcal{N}_2^{15r}({Y}_i+\frac{h}{2}[\![F(\tilde{Z}_1)]\!]_{2r})\Big\|_{L_q}+\frac{h}{2}\|[\![F(\tilde{Z}_1)]\!]_{2r}-F(\tilde{Z}_1)\|_{L_q} \nonumber
        \\&\leq\Big\|\mathcal{N}_2^{15r}(Y_{i,1})-\mathcal{N}_2^{15r}({Y}_i+\frac{h}{2}[\![F(\tilde{Z}_1)]\!]_{2r})\Big\|_{L_q}+\frac{h}{2}C_M\epsilon \label{eq:blubber}
    \end{align} 
    The second equality above holds because ${Y}_i+\frac{h}{2}[\![F(\tilde{Z}_1)]\!]_{2r}$ has rank at most $3r$, and therefore ${Y}_i+\frac{h}{2}[\![F(\tilde{Z}_1)]\!]_{2r}=\mathcal{N}_2^{15r}({Y}_i+\frac{h}{2}[\![F(\tilde{Z}_1)]\!]_{2r})$ holds almost surely. The last inequality follows from Lemma~\ref{lemma: low_rank_F}.
    With similar reasoning, we obtain the following bound for the first term in~\eqref{eq:blubber}:
    \begin{align*}
        \Big\|\mathcal{N}_2^{15r}(Y_{i,1})-\mathcal{N}_2^{15r}({Y}_i+\frac{h}{2}[\![F(\tilde{Z}_1)]\!]_{2r})\Big\|_{L_q} 
        =&\Big\|\mathcal{N}_2^{15r}(Y_{i,1})-({Y}_i+\frac{h}{2}[\![F(\tilde{Z}_1)]\!]_{2r})\Big\|_{L_q} 
        \\ \leq & \Big\|\mathcal{N}_2^{15r}(Y_{i,1})-({Y}_i+\frac{h}{2}F(\tilde{Z}_1))\Big\|_{L_q}+\frac{h}{2}C_M\epsilon.
    \end{align*}
    Using that $\hat{Z_1}=\tilde{Z_1}=Y_i$ holds almost surly, the law of total expectation, and Theorem~\ref{theorem: ny} give
    \begin{align*}
        \Big\|\mathcal{N}_2^{15r}(Y_{i,1})-({Y}_i+\frac{h}{2}F(\tilde{Z}_1))\Big\|_{L_q}
        =&\left(\mathbb{E}\left[\mathbb{E}(\|\mathcal{N}_2^{15r}(Y_{i,1})-(Y_{i,1})\|_F^q\vert \hat{Z}_1)\right]\right)^\frac{1}{q}
        \\\leq &C_{\mathcal{N}} \Big\|Y_{i,1}-[\![ Y_{i,1}]\!]_{15r}\Big\|_{L_q}
        \\ \leq &C_{\mathcal{N}}\Big\|{Y}_i+\frac{h}{2}F({Y}_i)-[\![ {Y}_i+\frac{h}{2}F({Y}_i)]\!]_{3r}\Big\|_{L_q}\\ \leq & C_{\mathcal{N}}\Big\|{Y}_i+\frac{h}{2}F({Y}_i)-({Y}_i+\frac{h}{2}[\![F({Y}_i)]\!]_{2r})\Big\|_{L_q} \leq  \frac{C_{\mathcal{N}}h}{2}C_M\epsilon.
    \end{align*}
    In summary, we have proven
    \begin{align}
    \label{eq: stage 2}
        \|\hat{Z_2}-\tilde{Z_2}\|_{L_q}&\leq \frac{(C_{\mathcal{N}}+2)C_M}{2}h\epsilon.
    \end{align} Note that we also have the following inequality  by Lemma \ref{lemma: low_rank_F}:
    \begin{equation}
    \label{eq: stage 2_low rank}
        \|\tilde{Z_2}-[\![\tilde{Z_2}]\!]_{3r}\|_{L_q}\leq \Big\|{Y}_i+\frac{h}{2}F(\tilde{Z}_1)-{Y}_i-\frac{h}{2}[\![F(\tilde{Z}_1)]\!]_{2r}\Big\|_{L_q}\leq\frac{h}{2}C_M\epsilon.
    \end{equation} 
    
For $j = 3$, we set $Y_{i,2}:= {Y}_i+\frac{h}{2}F(\hat{Z}_2)$ and apply an analogous reasoning:
      \begin{align*}
       \|\hat{Z_3}-\tilde{Z_3}\|_{L_q}&=\Big\|\mathcal{N}_3^{15r}(Y_{i,2})-{Y}_i-\frac{h}{2}F(\tilde{Z}_2)\Big\|_{L_q} \nonumber
        \\&\leq \Big\|\mathcal{N}_3^{15r}(Y_{i,2})-\mathcal{N}_3^{15r}({Y}_i+\frac{h}{2}[\![F([\![\tilde{Z}_2]\!]_{3r})]\!]_{6r})\Big\|_{L_q}\nonumber
        \\&\quad+\Big\|\mathcal{N}_3^{15r}({Y}_i+\frac{h}{2}[\![F([\![\tilde{Z}_2]\!]_{3r})]\!]_{6r})-{Y}_i-\frac{h}{2}F(\tilde{Z}_2)\Big\|_{L_q}\nonumber
               \\&= \Big\|\mathcal{N}_3^{15r}(Y_{i,2})-\mathcal{N}_3^{15r}({Y}_i+\frac{h}{2}[\![F([\![\tilde{Z}_2]\!]_{3r})]\!]_{6r})\Big\|_{L_q}+\frac{h}{2}\|
         F(\tilde{Z}_2)-[\![F([\![\tilde{Z}_2]\!]_{3r})]\!]_{6r}\|_{L_q}\nonumber
          \\&= \Big\|\mathcal{N}_3^{15r}(Y_{i,2})-({Y}_i+\frac{h}{2}[\![F([\![\tilde{Z}_2]\!]_{3r})]\!]_{6r})\Big\|_{L_q}+\frac{h}{2}\|
         F(\tilde{Z}_2)-[\![F([\![\tilde{Z}_2]\!]_{3r})]\!]_{6r}\|_{L_q}\nonumber
          \\&\leq  \Big\|\mathcal{N}_3^{15r}(Y_{i,2})-(Y_{i,2})\Big\|_{L_q}+\frac{h}{2}\|F(\hat{Z}_2)-[\![F([\![\tilde{Z}_2]\!]_{3r})]\!]_{6r}\|_{L_q}\\&\quad +\frac{h}{2}\|
         F(\tilde{Z}_2)-[\![F([\![\tilde{Z}_2]\!]_{3r})]\!]_{6r}\|_{L_q}.\nonumber
         \end{align*}
         The first term of the last inequality can be bounded using Lemma \ref{lemma: low-rank bound}:
         \begin{align*}
             \Big\|\mathcal{N}_3^{15r}(Y_{i,2})-(Y_{i,2})\Big\|_{L_q} \leq C_{\mathcal{N}}\Big\|Y_{i,2}-[\![ Y_{i,2}]\!]_{15r}\Big\|_{L_q}\nonumber
             \leq \frac{h}{2}C_{\mathcal{N}}\|F(\hat{Z}_2)-[\![F([\![\tilde{Z}_2]\!]_{3r})]\!]_{6r}\|_{L_q}.
         \end{align*}
Therefore, by \eqref{eq: stage 2} and \eqref{eq: stage 2_low rank},
\begin{align*}
       \|\hat{Z_3}-\tilde{Z_3}\|_{L_q} \leq &\frac{(C_{\mathcal{N}}+1)h}{2}\|F(\hat{Z}_2)-[\![F([\![\tilde{Z}_2]\!]_{3r})]\!]_{6r}\|_{L_q}+\frac{h}{2}\|
         F(\tilde{Z}_2)-[\![F([\![\tilde{Z}_2]\!]_{3r})]\!]_{6r}\|_{L_q}\\
         \leq & \frac{(C_{\mathcal{N}}+2)h}{2}(\|F(\hat{Z}_2)-F(\tilde{Z}_2)\|_{L_q}+\|F(\tilde{Z}_2)-F([\![\tilde{Z}_2]\!]_{3r})\|_{L_q}\\
         &+\|F([\![\tilde{Z}_2]\!]_{3r})-[\![F([\![\tilde{Z}_2]\!]_{3r})]\!]_{6r}\|_{L_q})
       \\ \leq &\frac{h}{2}(C_{\mathcal{N}}+2)(L\frac{(C_{\mathcal{N}}+2)C_M}{2}h\epsilon+L\frac{h}{2}C_M\epsilon+C_M\epsilon)
       \leq C_3h\epsilon.
       \end{align*}
       Also, as for $j = 2$, 
    \begin{equation*}
        \|\tilde{Z_3}-[\![\tilde{Z_3}]\!]_{7r}\|_{L_q}\leq \Big\|{Y}_i+\frac{h}{2}F(\tilde{Z}_2)-{Y}_i-\frac{h}{2}[\![F([\![\tilde{Z}_2]\!]_{3r})]\!]_{6r}\Big\|_{L_q}\leq\frac{h}{2}[C_M\epsilon+L\frac{h}{2}C_M\epsilon].
    \end{equation*}
   
Finally, for $j = 4$, we set $Y_{i,3}:= {Y}_i+{h}F(\hat{Z}_3)$ and obtain
    \begin{align*}
        \|\hat{Z_4}-\tilde{Z_4}\|_{L_q}
        &=\left\|\mathcal{N}_4^{15r}(Y_{i,3})-\hat{Y}_i-{h}F(\tilde{Z}_3)\right\|_{L_q}
        \\&\leq \left\|\mathcal{N}_4^{15r}(Y_{i,3})-\mathcal{N}_4^{15r}({Y}_i+h[\![F([\![\tilde{Z}_3]\!]_{7r})]\!]_{14r})\right\|_{L_q}
        \\&\quad+\left\|\mathcal{N}_4^{15r}({Y}_i+{h}[\![F([\![\tilde{Z}_3]\!]_{7r})]\!]_{14r})-{Y}_i-{h}F(\tilde{Z}_3)\right\|_{L_q}
         \\&\leq\left\|\mathcal{N}_4^{15r}(Y_{i,3})-({Y}_i+h[\![F([\![\tilde{Z}_3]\!]_{7r})]\!]_{14r})\right\|_{L_q}+h\left\|F(\tilde{Z}_3)-[\![F([\![\tilde{Z}_3]\!]_{7r})]\!]_{14r}\right\|_{L_q}
         \\&\leq\left\|\mathcal{N}_4^{15r}(Y_{i,3})-(Y_{i,3})\right\|_{L_q}+h\left\|F(\hat{Z}_3)-[\![F([\![\tilde{Z}_3]\!]_{7r})]\!]_{14r}\right\|_{L_q}\\&\quad +h\left\|F(\tilde{Z}_3)-[\![F([\![\tilde{Z}_3]\!]_{7r})]\!]_{14r}\right\|_{L_q}.
    \end{align*}
             The first term of the last inequality can once again be bounded using Lemma \ref{lemma: low-rank bound}:
\begin{align*}
             \left\|\mathcal{N}^{15r}(Y_{i,3})-(Y_{i,3})\right\|_{L_q} \leq C_{\mathcal{N}}\|Y_{i,3}-[\![ Y_{i,3}]\!]_{15r}\|_{L_q}
             \leq {h}C_{\mathcal{N}}\|F(\hat{Z}_3)-[\![F([\![\tilde{Z}_3]\!]_{7r})]\!]_{14r}\|_{L_q}.
         \end{align*}  
    Hence,
    \begin{align*}
         \|\hat{Z_4}-\tilde{Z_4}\|_{L_q} \leq & {h}(C_{\mathcal{N}}+1)\|F(\hat{Z}_3)-[\![F([\![\tilde{Z}_3]\!]_{7r})]\!]_{14r}\|_{L_q}+h\|F(\tilde{Z}_3)-[\![F([\![\tilde{Z}_3]\!]_{7r})]\!]_{14r}\|_F
         \\
         \leq& {h(C_{\mathcal{N}}+2)}(\|F(\hat{Z}_3)-F(\tilde{Z}_3)\|_{L_q}+\|F(\tilde{Z}_3)-F([\![\tilde{Z}_3]\!]_{7r})\|_F\\
         &+\|F([\![\tilde{Z}_3]\!]_{7r})-[\![F([\![\tilde{Z}_3]\!]_{7r})]\!]_{14r}\|_F)
       \\ \leq & {h}(C_{\mathcal{N}}+2)(LC_3h\epsilon+L\frac{h}{2}[C_M\epsilon+L\frac{h}{2}C_M\epsilon]+C_M\epsilon)
\leq C_4h\epsilon.
    \end{align*}
    
Collecting the obtained bounds for the stages and using 
Lipschitz continuity, we have
    \begin{align}
          \Big\|&{Y}_i+\frac{h}{6}(F(\hat{Z}_1)+2F(\hat{Z}_2)+2F(\hat{Z}_3)+F(\hat{Z}_4))\nonumber \\
          &-{Y}_i-\frac{h}{6}(F(\tilde{Z}_1)+2F(\tilde{Z}_2)+2F(\tilde{Z}_3)+F(\tilde{Z}_4))\Big\|_{L_q}\leq C_5h^2\epsilon. \label{eq: rk_local_lq}
    \end{align} Recall that one step of the classic RK 4 method has error $\mathcal{O}(h^5)$. With the true solution $\Phi^{h}_F({Y}_{i})$ satisfying the low-rank approximability assumption~\eqref{eq:low-rank_assumption_adaptive}, we find that one step of RK 4 satisfies
    \begin{align}
        \Big\|&{Y}_i+\frac{h}{6}(F(\tilde{Z}_1)+2F(\tilde{Z}_2)+2F(\tilde{Z}_3)+F(\tilde{Z}_4)) \nonumber \\
        &-[\![{Y}_i+\frac{h}{6}(F(\tilde{Z}_1)+2F(\tilde{Z}_2)+2F(\tilde{Z}_3)+F(\tilde{Z}_4))]\!]_{r}\Big\|_F\leq C_6(h^5+h\epsilon). \label{eq: rk_local_low_rank}
        \end{align}
Using Theorem~\ref{theorem: ny} and the inequalities~\eqref{eq: rk_local_lq}, \eqref{eq: rk_local_low_rank}, we obtain \begin{align*}
        &\Big\|{Y}_{i+1}-{Y}_i-\frac{h}{6}(F(\hat{Z}_1)+2F(\hat{Z}_2)+2F(\hat{Z}_3)+F(\hat{Z}_4))\Big\|_{L_q}\\
        \leq & C_\mathcal{N}\Big\|{Y}_i+\frac{h}{6}(F(\hat{Z}_1)+2F(\hat{Z}_2)+2F(\hat{Z}_3)+F(\hat{Z}_4))\\
        & \quad\ \ -[\![{Y}_i+\frac{h}{6}(F(\hat{Z}_1)+2F(\hat{Z}_2)+2F(\hat{Z}_3)+F(\hat{Z}_4))]\!]_{r}\Big\|_{L_q}\\
        \leq & C_\mathcal{N}\Big\|{Y}_i+\frac{h}{6}(F(\hat{Z}_1)+2F(\hat{Z}_2)+2F(\hat{Z}_3)+F(\hat{Z}_4))\\
        & \quad\ \ -[\![{Y}_i+\frac{h}{6}(F(\tilde{Z}_1)+2F(\tilde{Z}_2)+2F(\tilde{Z}_3)+F(\tilde{Z}_4))]\!]_{r}\Big\|_{L_q}\\& \leq C_{\mathcal{N}}(C_5h^2\epsilon+C_6(h^5+h\epsilon)).
    \end{align*} The result of the lemma is concluded from the following inequality:
    \begin{align*}
        \|\Phi^{h}_F({Y}_{i})-{Y}_{i+1}\|_{L_q}&\leq \|\Phi^{h}_F({Y}_{i})-{Y}_i-\frac{h}{6}(F(\tilde{Z}_1)+2F(\tilde{Z}_2)+2F(\tilde{Z}_3)+F(\tilde{Z}_4))\|_{L_q}\\&\quad +C_5h^2\epsilon+C_{\mathcal{N}}(C_5h^2\epsilon+C_6(h^5+h\epsilon)) \leq \mathcal{O}(h^5+h\epsilon).
    \end{align*}
\end{proof}
Finally, the following theorem establishes order $4$ with respect to $h$ of the modified randomized low-rank RK 4 method~\eqref{eq: modified scheme}. This provides some theoretical explanation for the convergence order $4$ we observe for the randomized low-rank RK 4 method (without intermediate rank increases).
\begin{theorem} \label{theorem:rk4}
Suppose that the assumptions stated in Section~\ref{section: assumptions}
and~\eqref{eq:low-rank_assumption_adaptive} hold. 
Under the assumption \eqref{eq:low-rank_assumption_adaptive} and the assumptions stated in section \ref{section: assumptions}. The the global error of the scheme~\eqref{eq: modified scheme} with independent standard Gaussian matrices and oversampling parameters $p,\ell\geq 4$ satisfies for $q=\min\{p,\ell\}$ the bound
    \begin{equation*}
         \|{Y}_{N}-A(Nh)\|_{L_q}\leq C(\epsilon+h^4+\delta),
    \end{equation*}on the finite time-interval $0 \leq nh \leq T$ and for every $0 < h \leq h_0$. The constant $C$ depends only on $L,T,h_0$ $C_M$ and $C_{\mathcal{N}}$.
    In particular, for any $\eta \geq 1$, it holds for fixed $h$ and $N$ that
\begin{equation*}
    \Pr\{ \|Y_N-A(Nh)\|_F\geq C\eta(\epsilon+h^4+\delta)\}\leq\frac{1}{\eta^q}.
\end{equation*}
\end{theorem} 
\begin{proof}
    By Lemma \ref{lemma: local error of Yhat} we have 
    \begin{equation*}
        \Big(\mathbb{E}(\mathbb{E}[\|\Phi^{h}_F({Y}_{i})-{Y}_{i+1}\|_F^q|Y_i])\Big)^{1/q}\leq \Big(\mathbb{E}[C^q(h^5+h\epsilon)^q]\Big)^{1/q}= C(h^5+h\epsilon).
    \end{equation*} Substituting this bound into~\eqref{eq: bound of local error} proves the first statement of the theorem. The second statement is proved by Markov's inequality, as in the proof of Corollary \ref{corollary: idealised_method_prob}.
\end{proof}

\section{Numerical Experiments}
\label{sect: Numerical Experiments}

In the following numerical experiments, we verify the accuracy of randomized low-rank RK methods, Algorithm \ref{alg:rrk-s}, which will be denoted as Rand RK. In particular, we consider Rand RK1 (Euler), Rand RK2 and Rand RK4, based upon
the Euler method, Heun's method and the classical RK 4 method, respectively. For convenience, we recall the corresponding coefficients:
\begin{itemize}
    \item Rand RK1 (Euler): $b_1=1$
    \item Rand RK2: $a_{21}=1$, $b_1=b_2=\frac{1}{2},$
    \item Rand RK4: $a_{21}=a_{32}=\frac{1}{2}$, $a_{43}=1,$ $b_1=b_4=\frac{1}{6}$, $b_2=b_3=\frac{1}{3}$.
\end{itemize}
We have implemented the generalized Nyst\"{o}m method following~\cite{tropp2017fixed}; we always use Gaussian random matrices and the oversampling parameters $p = \ell = \max\{2,0.1r\}$ for sketching. Although the generalized Nystr\"{o}m method is usually numerically stable~\cite{nakatsukasa2020fast}, it may exhibit instabilities, especially when $\Psi^TZ\Omega$ is numerically rank deficient.  We have never observed such an instability in our numerical experiments, but let us point out that a numerically safer
variant with regularization is described in \cite{nakatsukasa2020fast}. Note that our implementation of Rand RK4 does \emph{not} use oversampling in the intermediate stage, that is, we have implemented Algorithm~\ref{alg:rrk-s} with the coefficients of RK 4 instead of \eqref{eq: modified scheme}.

The error is measured by the Frobenius norm between the reference solution and the approximation. The reference solution is obtained by solving the full matrix differential equation with MATLAB \verb+ODE45+ using the tolerances \verb+{’RelTol’, 1e-10, ’AbsTol’, 1e-10}+.  Because our methods involve randomization, we report the mean approximation error as well as the spread between the largest and smallest errors for 10 independent random trials (indicated by lower/upper horizontal lines in the graph).

In the first two experiments, we also report the error of our implementation of the projected RK $s$ method~\cite{kieri2019projection} for $s=1,2,4$  as well as the projector
splitting integrator~\cite{Lubich2014projector} with the sub-steps computed by MATLAB's \verb+ODE45+ using the same tolerance parameters as above. The initial value for these methods is obtained by applying the truncated SVD to the initial matrix $A_0$.

All experiments have been performed in Matlab (version 2023a) on a Macbook Pro with an Apple M1 Pro processor. The code used to perform the experiments and produce the figures can be found at \url{https://github.com/hysanlam/rand_RK}.

\subsection{Lyapunov matrix differential equation} \label{sec:lyapunovexample}

As a first simple experiment, we approximate the solution of a Lyapunov matrix differential equation~\cite[Section 6.1]{Uschmajew2020Geometric}, which takes the form
\begin{equation*}
    \dot{A}(t)=LA(t)+A(t)L+\alpha \frac{C}{\|C\|_F}, \quad A(0)=A_0,
\end{equation*}
 where $A(t)\in \mathbb{R}^{n\times n}$, $L\in \mathbb{R}^{n\times n}$, $C\in \mathbb{R}^{n\times n}$, $\alpha\geq 0$ and $t\in [0, T]$.
 We set $n=128$ and use the symmetric matrix $L=\text{diag}(1,-2,1)\in \mathbb{R}^{128 \times 128 }$. For setting the entries of the source term $C$ and initial matrix $A_0$, we follow a construction similar to the one used in~\cite[section 5.1]{ceruti2024robust}:
\begin{equation*}
   C_{ij}=\sum^{11}_{k=1}10^{-(k-1)}\cdot e^{-k(x_i^2+y_j^2)},
\end{equation*} $$(A_0)_{ij}=\sum^{20}_{k=1}b_k\cdot \sin(kx_i) \sin(ky_j),\text{ with }b_k=\begin{cases}
1, \text{ if }k=1\\
5e^{-(7+0.5(k-2))},\text{ if }k>1,\\
\end{cases}$$ where $(x_i,y_j)$, for $i,j=1,\ldots 128$, are uniform discretization points of the square $[-\pi,\pi]\times [-\pi,\pi]$. The final time is set to $T=1$ and first consider $\alpha=1$. Figure \ref{Figure: Lyapunov RK} displays the singular values of the reference solution at $t=T$, and the error of the approximation obtained by Rand Euler and Rand RK4 with different ranks. We can see that Rand Euler achieves first-order convergence in time, while Rand RK4 achieves fourth-order convergence in time until it reaches the level of low rank approximation error. Moreover, both methods demonstrate robust behavior despite randomness; among these 10 trials, the maximum error is only at most three times the empirical mean.\begin{figure}[H]
 \makebox[\textwidth][c]{\includegraphics[scale=0.3]{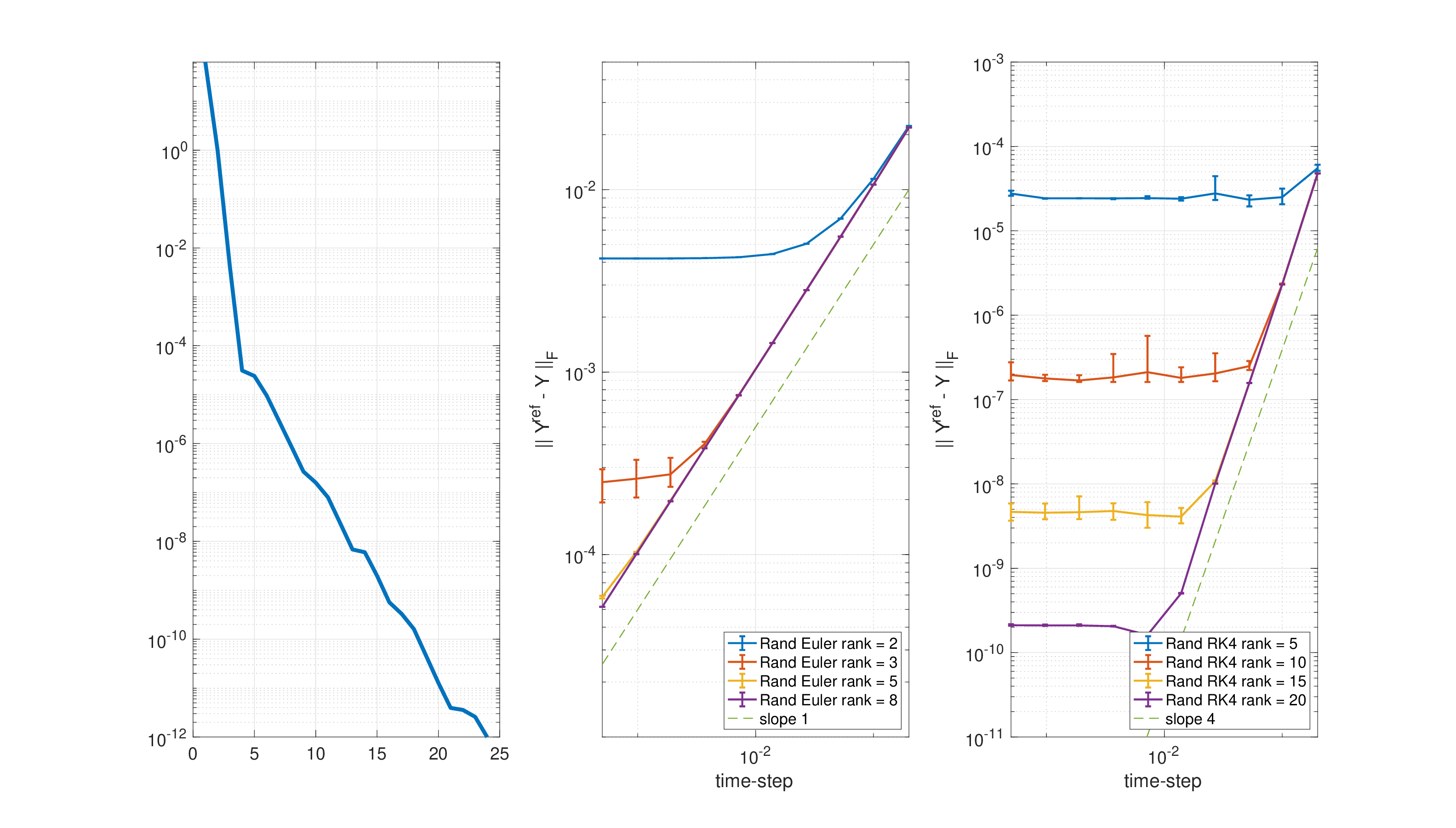}}
 \caption{Lyapunov matrix differential equation with $\alpha=1$. The singular values of the reference solution at time $T = 1$ together with the approximation
errors of the numerical approximation obtained via the Rand Euler and Rand RK4 for different ranks and time-step sizes.}
 \label{Figure: Lyapunov RK}
\end{figure}
We now fix the rank to $r=10$ and compare the accuracy of Rand Euler, Rand RK2, and Rand RK4 with other methods, including PRK1, PRK2, PRK4, and projector splitting for $\alpha=10^{-5}$ as well as $\alpha=1$. In Table \ref{Table: proj_error_1}, we report the average and maximum values of the tangential projection error $\|F(Y_i)-P_r{(Y_i)}F(Y_i)\|_F$, where $Y_i$ is the approximation computed at the $i$th time step by PRK 2 with $h=5\times 10^{-3}$. For $\alpha=1$ this error is much larger than $\alpha=10^{-5}$. However, when the tangential projection error is large, the accuracy and convergence behavior of PRK and projector splitting are not guaranteed. This is indeed observed in Figure~\ref{fig: Lyapunov 1}. For $\alpha=10^{-5}$, PRK 1, PRK 2, and projector splitting exhibit the expetected order of convergence. PRK 4 seems to even exhibit fourth-order convergence initially but this quickly deteriorates for smaller $h$. When $\alpha=1$, all these methods only show first-order convergence. On the other hand, the randomized methods remain robust: Rand Euler, Rand RK2, and Rand RK4 exihibit first, second, and fourth-order convergence, respectively, for both choices of $\alpha$.
\begin{table}[H]
\centering
\begin{tabular}{|l|l|l|}
\hline
 $\alpha$ &$10^{-5}$  & 1\\ \hline\hline 
average $\|F(Y_i)-P_r{(Y_i)}F(Y_i)\|_F$ &$1.3553\times 10^{-7}$  &  $4.0144\times 10^{-4}$ \\ \hline
max $\|F(Y_i)-P_r{(Y_i)}F(Y_i)\|_F$&$1 \times 10^{-5}$   &  $0.7932$ \\ \hline
\end{tabular}
 \caption{Lyapunov matrix differential equation with $\alpha=10^{-5}$ and $\alpha=1$. Average and maximum $\|F(Y_i)-P_r{(Y_i)}F(Y_i)\|_F$ for the approximation $Y_i$ at the $i$th time step of PRK 2 with $h=5\times 10^{-3}$. }
 \label{Table: proj_error_1}
\end{table}
\begin{figure}[H]
 \makebox[\textwidth][c]{\includegraphics[scale=0.32]{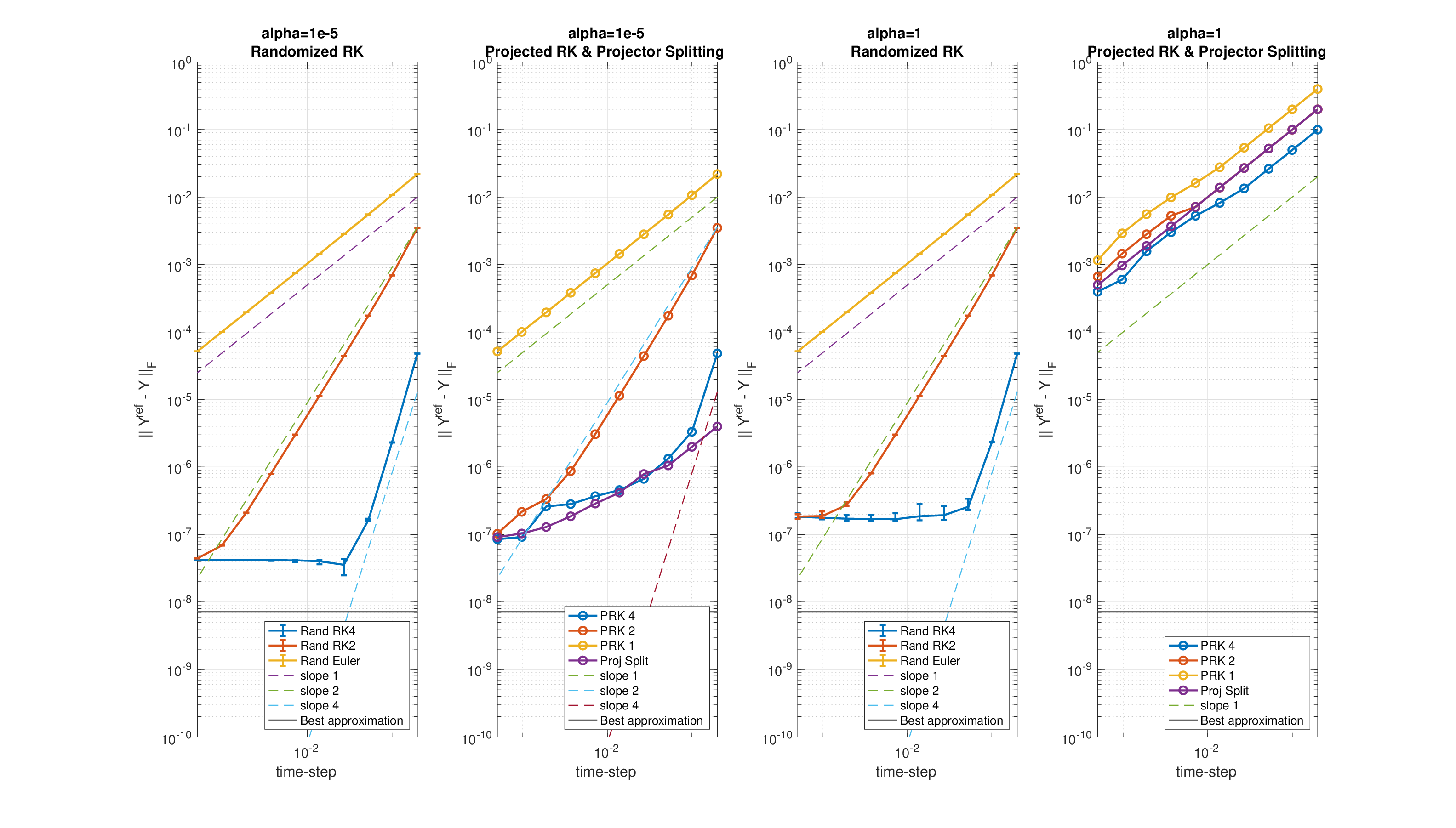}}
 \caption{Lyapunov matrix differential equation with $\alpha=10^{-5}$ and $\alpha=1$. Comparison of absolute approximation errors for
different low-rank integrators using rank $r = 10$.}
  \label{fig: Lyapunov 1}
\end{figure}

\subsubsection{Speeding up by constant random matrices}
\label{section 4.1.1}

As indicated in Section~\ref{sec:compprk}, the cost of computing the sketches $\hat{K}_{jq},\tilde{K}_{jq}$ in Algorithm~\ref{alg:rrk-s}
can be divided by nearly a factor $s$ when choosing the same random matrices across the stages in one time step. That is, we draw two independent Gaussian matrices $\Omega_1$ and $\Psi_1$, and set $\Omega_1=\Omega_2=\cdots=\Omega_{s+1}$, $\Psi_1=\Psi_2=\cdots=\Psi_{s+1}$. In this experiment, we ran 10 random trials to solve the differential Lyapunov equation for $\alpha = 1$ and compare the described constant choice with the standard (independent) choice of random matrices in Algorithm \ref{alg:rrk-s}. We tested Rand RK2 and Rand RK4 with different ranks and plotted the obtained errors in Figure~\ref{fig: Lyapunov 3}. For this example, we observe comparable performance when using either different or the same random matrices across the stages in a time step. Although this modification is computationally attractive and likely the preferred way to run Algorithm \ref{alg:rrk-s}, we are unable to establish an error bound due to the lack of stochastic independence of the stages.
\begin{figure}[H]
 \makebox[\textwidth][c]{\includegraphics[scale=0.32]{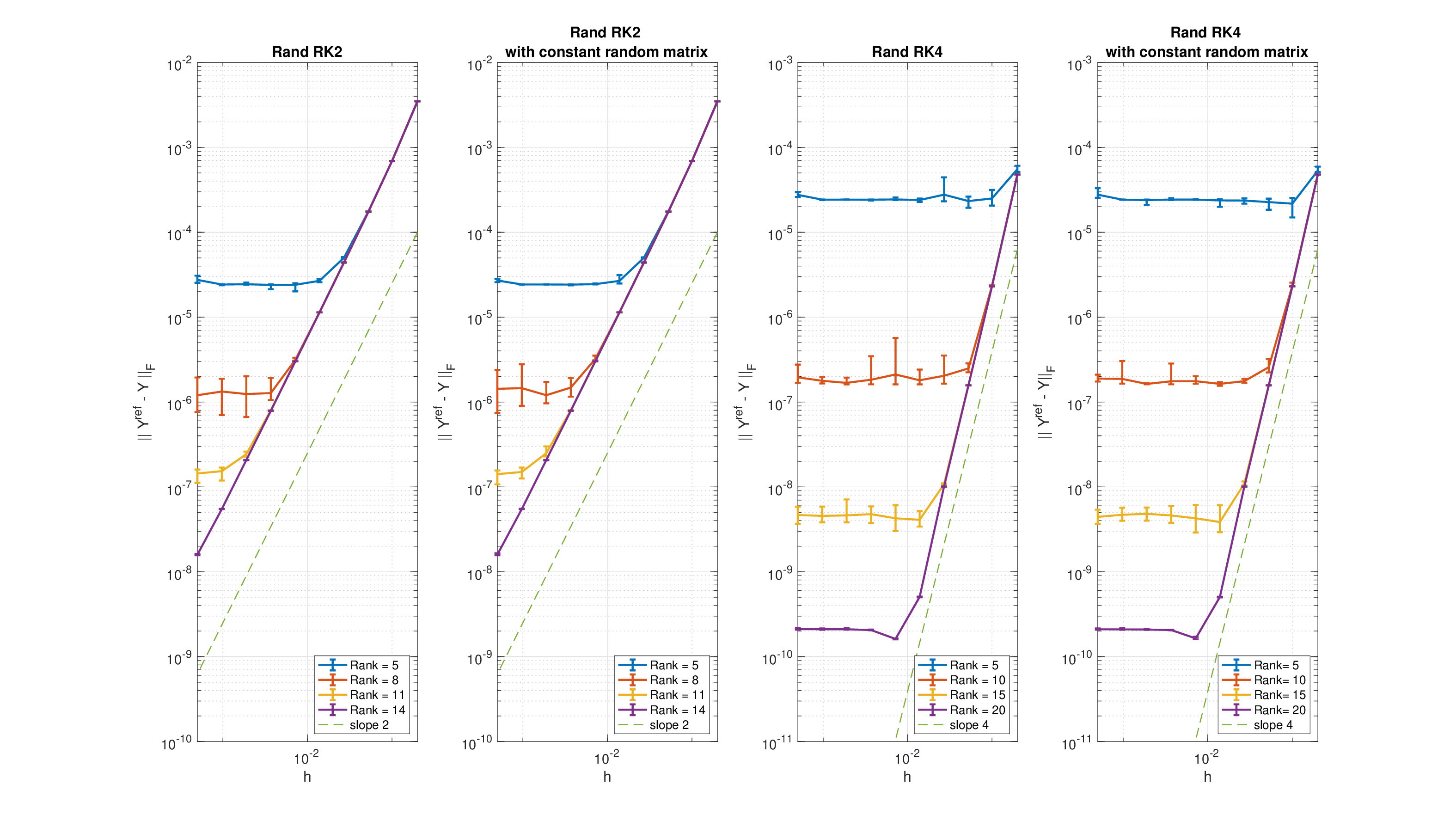}}
\caption{Lyapunov matrix differential equation with $\alpha=1$. Comparison of absolute approximation errors for Rand RK2, Rand RK2 with different (independent) random matrices across stages vs. Rand RK2 and Rand RK4 with the same random matrices across stages.}
\label{fig: Lyapunov 3}
\end{figure}

\subsection{Non-linear Schr\"{o}dinger equation}

We now consider the non-linear Schr\"{o}dinger equation from~\cite[Section 5.3]{kieri2019projection}, where $A:[0,T]\rightarrow \mathbb{C}^{n\times n}$ evolves according to\begin{equation}
     \dot{A}(t)=i[\frac{1}{2}(BA+AB)+\alpha|A|^2A].
\end{equation} The cubic nonlinearity $|A|^2A$ is taken element-wise and $B=\mathrm{diag}(1,0,1)$. We choose $n=100$, $T=5$ and the initial data  $$(A_0)_{ij}=\exp\Big(-\frac{(i-60)^2}{100}-\frac{(j-50)^2}{100}\Big)+\exp\Big(-\frac{(i-50)^2}{100}-\frac{(j-40)^2}{100}\Big).$$ In this example, we aim at computing approximations of rank up to $30$. To ensure that $A_0$ has rank at least $30$ (making sure it satisfies the assumptions of PRK), we perturb $A_0$ by taking the full SVD of $A_0$ and setting the singular values $3,4,\ldots, 32$ to $10^{-9}$. 

First, we set $\alpha=0.3$. The singular values of the reference solution at $t=T$, as well as the approximation errors by Rand Euler and Rand RK4 with different ranks are shown in Figure \ref{Figure: Non-linear 1}. We again observe that Rand Euler achieves first-order convergence in time, while Rand RK4 achieves fourth-order convergence in time until it reaches the level of low-rank approximation error. In this example, the maximum error of both methods is also very close to the empirical mean, deviating by less than twice the mean.
\begin{figure}[H]
 \makebox[\textwidth][c]{ \includegraphics[scale=0.3]{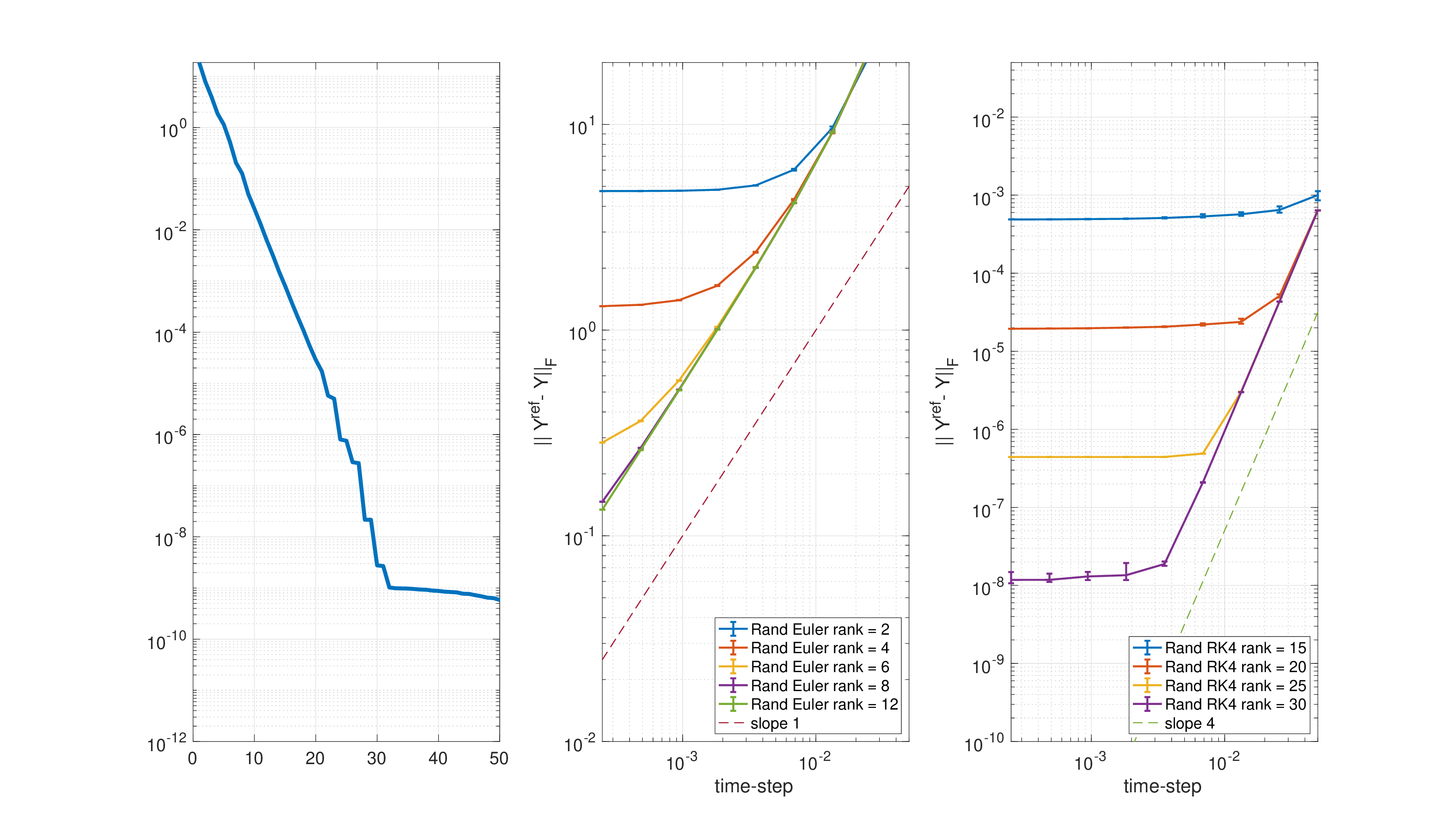}}
  \caption{Non-linear Schr\"{o}dinger equation with $\alpha=0.3$. Singular values of the reference solution at time $T = 5$ together with the approximation
errors of the numerical approximation obtained by Rand Euler and Rand RK4 for different ranks and time-step sizes. }
  \label{Figure: Non-linear 1}
\end{figure}

{Now we fix the rank to $r=30$ and compare Rand RK with PRK and projector splitting for $\alpha=3\times 10^{-4}$ and $ \alpha=0.3$ in Figure \ref{Figure: Non-linear 2}. For small  $\alpha=3\times 10^{-4}$, we observe that PRK 1 and PRK 2 exhibit the correct order of convergence. Again, one observes that the initially visible fourth-order convergence  of PRK 4 quickly deteriorates as  $h$ decreases. For $\alpha=0.3$, we see that the order of PRK 2 decreases to 1 when $h$ is small, and PRK 4 only shows first-order convergence as well. Again, this is likely due to the large tangential projection error $\|P_r{(Y_i)}F(Y_i)-F(Y_i)\|_F$; see Table \ref{table: non-linear}. On the other hand, all the Randomized RK methods exhibit robust convergence of the expected order. Surprisingly and for reasons unclear to us, the projector splitting method provides very accurate results for this example. }
\begin{table}[H]
\centering
\begin{tabular}{|l|l|l|}
\hline
 $\alpha$ &3e-4  & 3e-1\\ \hline\hline 
average $\|F(Y_i)-P_r{(Y_i)}F(Y_i)\|_F$ &$3.0145\times 10^{-8}$   &  $2.9315\times 10^{-5}$ \\ \hline
max $\|F(Y_i)-P_r{(Y_i)}F(Y_i)\|_F$&  $5.8059\times 10^{-4}$& $0.5806$\\ \hline
\end{tabular}
\caption{Average and maximum $\|F(Y_i)-P_r{(Y_i)}F(Y_i)\|_F$ for the approximation $Y_i$ at the $i$th time step of PRK 2 with $h= 2.5\times 10^{-4}$.} 
\label{table: non-linear}
\end{table}
\begin{figure}[ht]
 \makebox[\textwidth][c]{\includegraphics[scale=0.32]{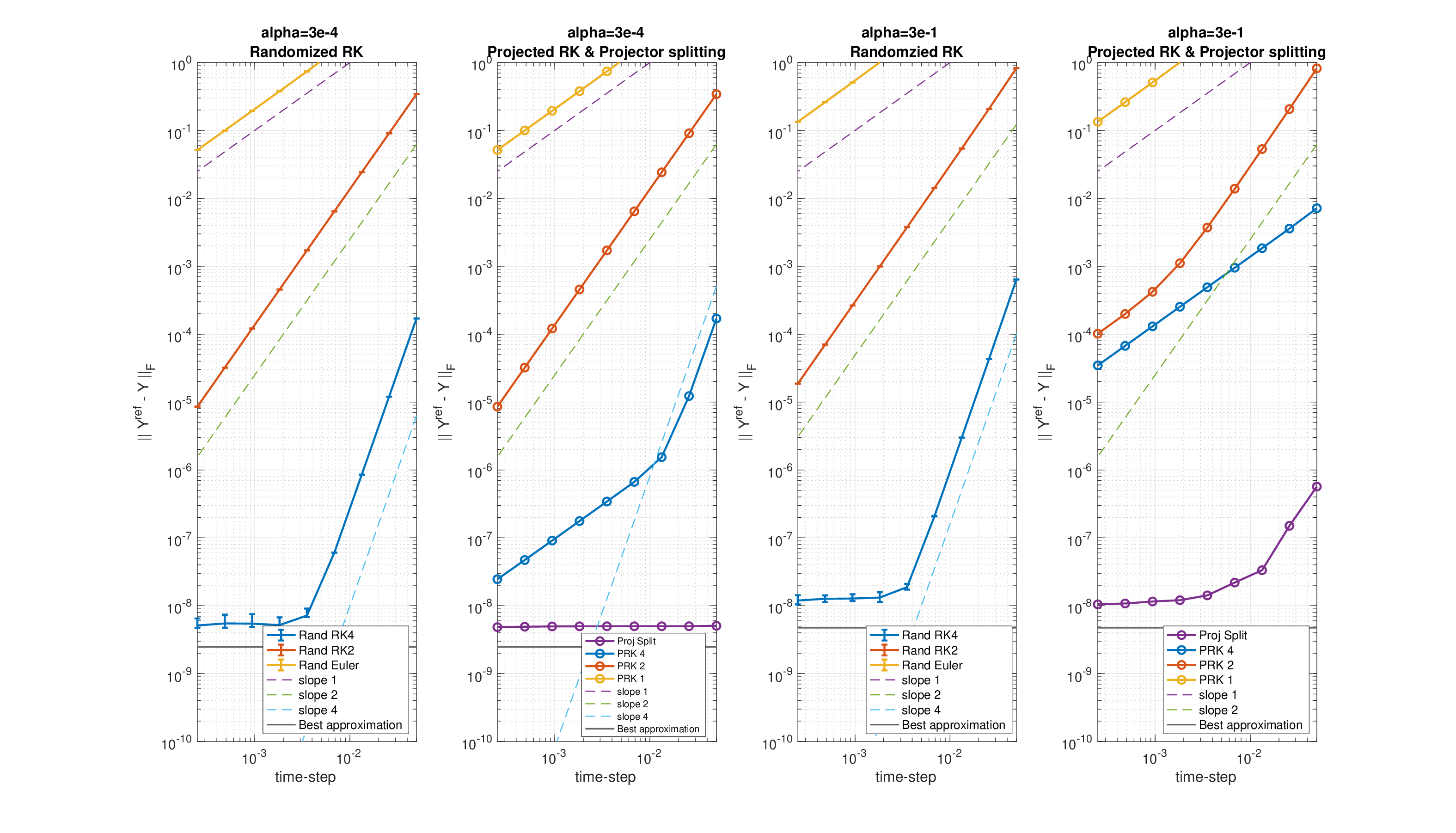}}
  \caption{Non-linear Schr\"{o}dinger equation with $\alpha=3\times 10^{-4}$ and $\alpha=3\times 10^{-1}$. Comparison of absolute approximation errors measured in Frobenius norm for
different integrators for rank-30. }
  \label{Figure: Non-linear 2}
\end{figure}

\subsection{Discrete Schr\"{o}dinger equation in imaginary time}

In this example, we aim at approximating the solution of the discrete Schr\"{o}dinger equation in
imaginary time from~\cite{ceruti2022unconventional}:
\begin{equation}\label{eq:discschroedinger}
    \dot{A}(t) =-H[A(t)],\quad A(0)=A_0,\quad t\in[0,0.5]
\end{equation} where
\begin{equation*}
    H[A(t)]=-\frac{1}{2}(DA(t)+A(t) D)+V_{\cos} A(t) V_{\cos}\in \mathbb{R}^{n\times n},
\end{equation*}
$D=\text{diag}(-1,2,-1)$ is the discrete 1D Laplace, and $V_{\cos}$ is the diagonal matrix with diagonal entries
$1-\cos(2j\pi/n)$ for $j=-n/2,\ldots,n/2-1$.
We choose $n=512$ and an initial value $A_0$ that is randomly generated with prescribed singular values $10^{-i}$, $i=1,\ldots,512$. In Figure \ref{fig: Discrete Schrodinger}, we first plot the singular values of the reference solution at $T=0.5$. Next, we plot the error of Rand Euler, Rand RK2, and Rand RK4 with a rank of 40. The final plot shows the error of Rand RK4 with various ranks. Once again, we observe the expected order of convergence until it reaches the level of low-rank approximation error.
\begin{figure}[h]
 \makebox[\textwidth][c]{ \includegraphics[scale=0.3]{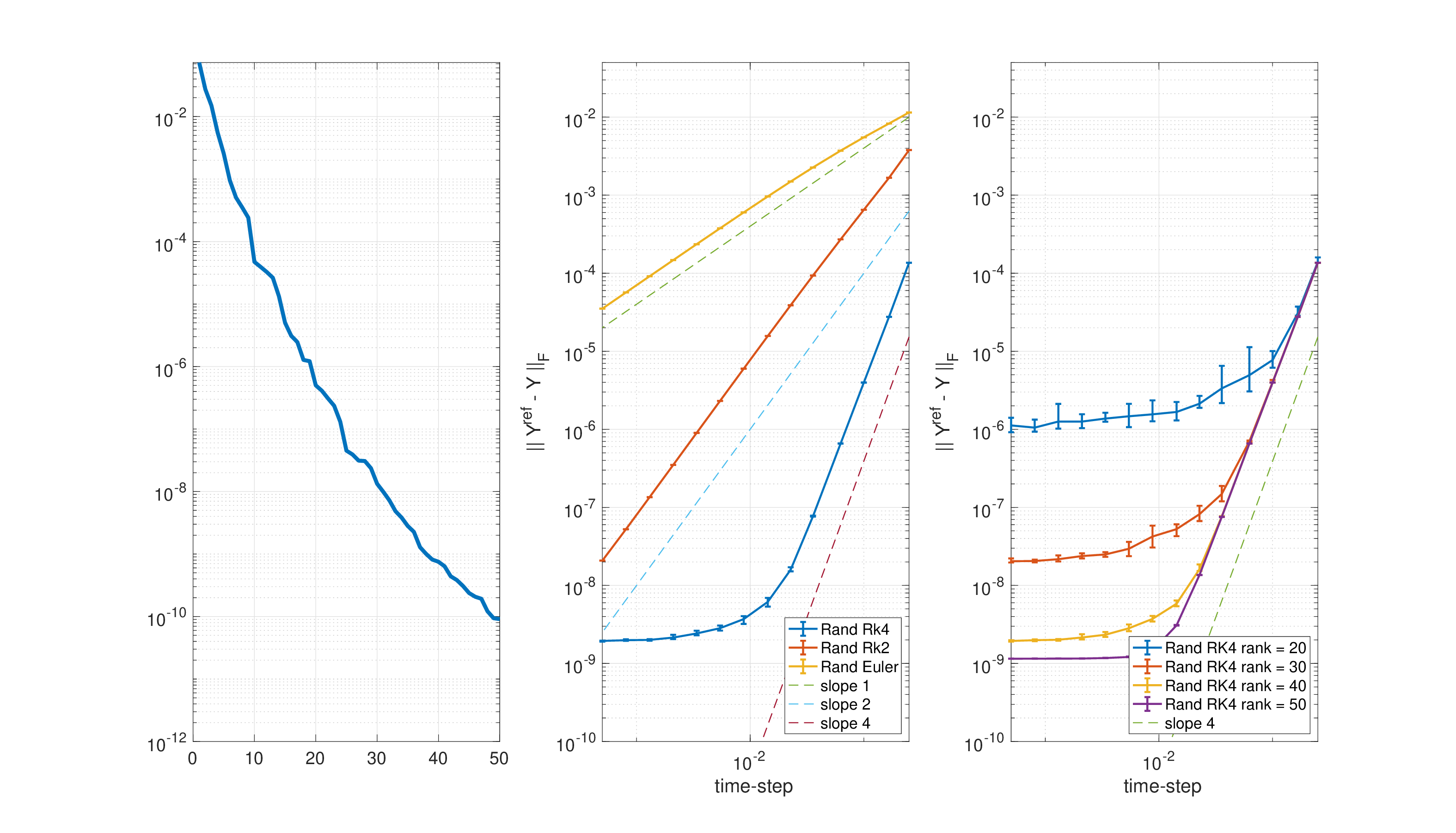}}
 \caption{Discrete Schr\"{o}dinger equation in imaginary time. The singular values of the reference solution at time $T = 0.5$, and rank-$40$ approximation errors for Rand Euler, Rand RK2 and Rand RK4. Also, the approximation errors for
different ranks using Rand RK4.}
\label{fig: Discrete Schrodinger}
\end{figure}

\subsection{Allen-Cahn equation}

Following~\cite[Section 5.3]{carrel2023projected}, we consider the matrix differential equation arising from discretizeing the Allen-Cahn equation via finite differences:
\begin{equation*}
    \dot{A}=\epsilon(LA+AL)+A-A^3,
\end{equation*} with initial data
\begin{equation*}
    (A_0)_{ij}=\frac{[e^{-\tan^2(x_i)}+e^{-\tan^2(y_j)}]\sin(x_i)\sin(y_j)}{1+e^{|csc(-x_i/2)|}+e^{|csc(-y_i/2)|}},
\end{equation*} where $A^3$ is to be understood element-wise and $(x_i,y_j)\in [0,2\pi]^2$, with $i,j=1,\ldots, 256$, are uniform discretization points. The matrix $L$ is the one-dimensional finite-difference stencil, $\epsilon=0.01$ and the time interval is $[0,10]$. For this example, we apply Rand RK4 with rank 2 and time step size $h=10^{-3}$. We plot the contours of the results at $t=1,3,5,7,10$ in Figure \ref{fig: allen-cahn}. Additionally, we calculate the difference between the reference solution and normalize the difference by the Frobenius norm of the reference solution. We observe that Rand RK 4 accurately captures the contours despite using a very low rank.
\begin{figure}[h]
\centering
\makebox[\textwidth][c]{ \includegraphics[scale=0.5]{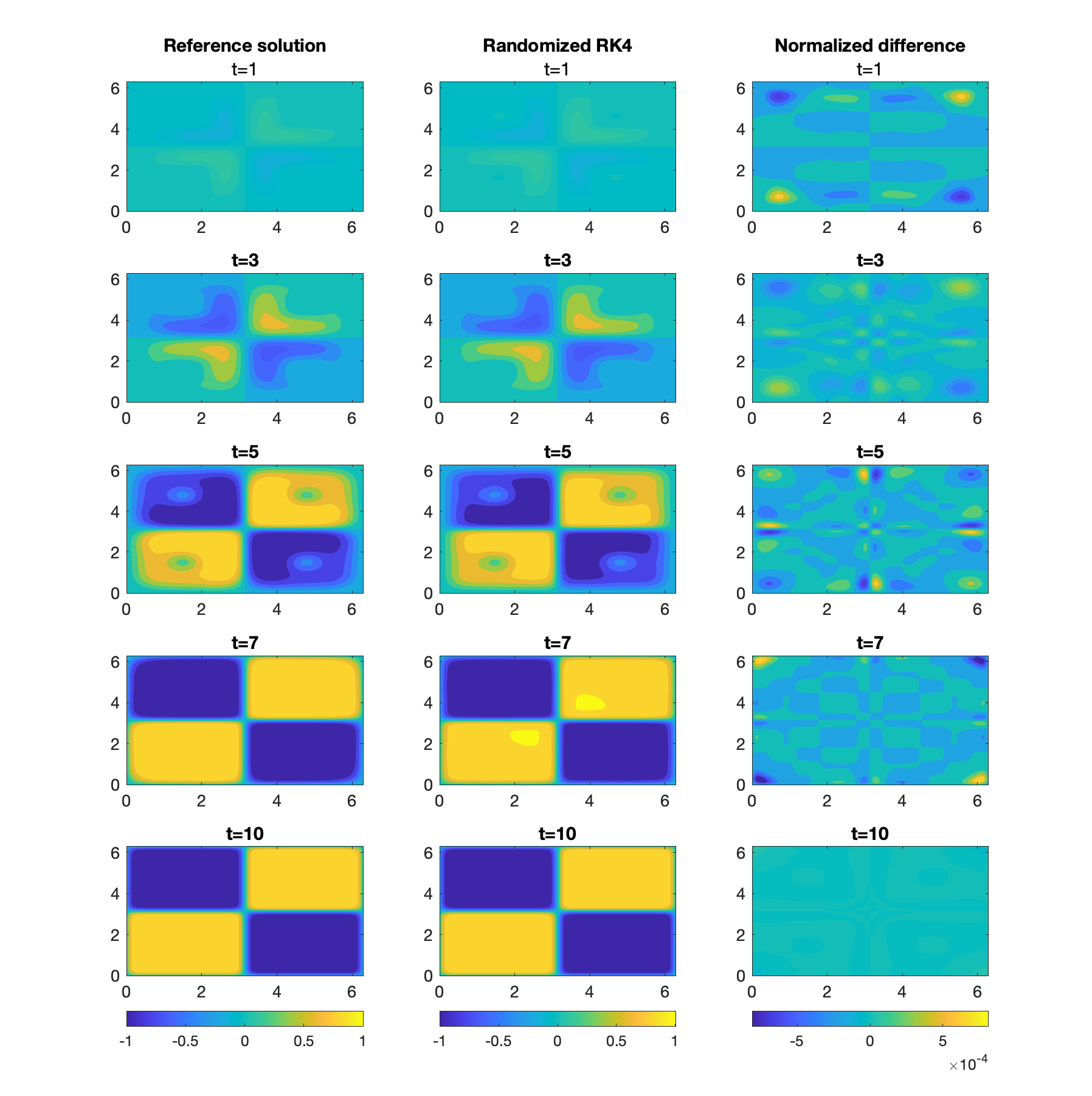}}
  \caption{Allen-Cahn equation.The contour of the reference solution and the rank-2 approximation obtained by Rand RK4 with $h = 10^{-3}$.}
  \label{fig: allen-cahn}
\end{figure}

\section{Conclusions}

In this work, we have proposed randomized low-rank Runge-Kutta methods. The analysis and numerical experiments clearly demonstrate the great potential of these methods to
constitute an attractive alternative to existing dynamical low-rank methods. To fully realize this potential, further work is needed. On the practical side, rank adaptivity, parallelization, the preservation of quantities, and the combination with splitting/exponential integrators belong to the aspects that need to be studied in order to match the progress on dynamical low-rank methods achieved during the last decade. Also, the use of structured random matrices for efficiently addressing matrix differential equations with nonlinearities merits exploration. On the theoretical side, our analysis has raised a number of open questions. In particular, it would be important to develop a more direct approach for establishing the full order of randomized low-rank Runge-Kutta methods.

\appendix 
\section{Appendix: Proof of Theorem \ref{Thm: stage order}}

The proof of Theorem~\ref{Thm: stage order} follows from modifying the proof of~\cite[Theorem 6]{kieri2019projection} appropriately. We first provide an bound that relates the stages of the standard RK method~\eqref{eq: standard RK} with the ones of the randomized RK method~\eqref{rand_rk}.
\begin{theorem}
With the assumptions stated in Theorem \ref{Thm: idealised lq}, suppose that  $\gamma_1\leq \gamma_2\leq \cdots \leq \gamma_s$ are the stage orders of the explicit RK method~\eqref{eq: standard RK} and denote $\tilde{\gamma}_j=\min(\gamma_j,\gamma_2+1)$. Then the differences between the stages $\tilde{Z}_j$ of the standard RK method~\eqref{eq: standard RK} with $A_i = Y_i$
and the stages $Z_j$  of the randomized RK method~\eqref{rand_rk} are bounded by
    \begin{equation} \label{eq:firsteqn}
      \|Z_j-\tilde{Z}_j\|_{L_q}\leq \begin{cases}
           0\quad &\text{ if }j=1,2,\\
           C(h^2\epsilon+h^{\gamma_2+2})&\text{ if }j=3,4,\ldots, s,
       \end{cases}
    \end{equation}
     \begin{equation}
     \label{eq: F_diff}
       \|F(\mathcal{N}_j(Z_j))-F(\tilde{Z}_j)\|_{L_q}\leq \begin{cases}
           0\quad &\text{ if }j=1,\\
           C(h\epsilon+h^{\tilde{\gamma}_j+1})&\text{ if }j=2,3,\ldots, s,
       \end{cases}
    \end{equation}
    for any $0\leq h\leq h_0$, where $q=\min\{p,\ell\}$ and the constant $C$ depends only on $L$, $T$, $C_\mathcal{N}$, $h_0$, $s$ and $\max_{ij}|a_{ij}|$.
\end{theorem}
\begin{proof}
    For $j=1$, we have $Y_i=Z_1=\tilde{Z}_1$ and, therefore, 
    \begin{equation*}
        F(\mathcal{N}_1(Z_1))=F(Z_1)=F(\tilde{Z}_1),
    \end{equation*} holds almost surely, because $Y_i$ has rank at most $r$.  For $j\geq 2$, we proceed by induction and assume that the statement of the theorem holds up to $j-1$. Using the induction hypothesis for~\eqref{eq: F_diff}, we obtain that
    \begin{align*}
       \|Z_j-\tilde{Z}_j\|_{L_q} & \leq h\sum^{j-1}_{l=1}|a_{jl}|\cdot \|F(\mathcal{N}_l(Z_l))-F(\tilde{Z}_l)\|_{L_q} \\
       &\leq \begin{cases}
            0\quad &\text{ if }j=2,\\
            C_ACh(sh\epsilon+h^{\tilde{\gamma}_2+1}+\ldots +h^{\tilde{\gamma}_{j-1}+1})&\text{ if }j=3,\ldots , s.
        \end{cases}
    \end{align*}
    Using the assumption on the ordering of the stage orders and $\tilde{\gamma}_2=\gamma_2$, it follows that
    \begin{equation}
    \label{eq: local_error_1}
          \|Z_j-\tilde{Z}_j\|_{L_q}\leq\begin{cases}
            0\quad &\text{ if }j=1,2,\\
            C_Zh(h\epsilon+h^{{\gamma}_2+1})&\text{ if }j=3,\ldots , s,
        \end{cases}
    \end{equation}
    which shows~\eqref{eq:firsteqn}.
    To establish~\eqref{eq: F_diff}, we note that $Z_j$ is independent from $\Omega_j$ and $\Psi_j$, which allows us to use the law of total expectation and Theorem~\ref{theorem: ny} to conclude that \begin{align*}
        \|F(\mathcal{N}_j(Z_j))-F(\tilde{Z_j})\|_{L_q}&\leq L\|\mathcal{N}_j(Z_j)-\tilde{Z_j}\|_{L_q}
      \\&\leq L\left(\big(\mathbb{E}\{\mathbb{E}[\|\mathcal{N}_j(Z_j)-{Z_j}\|^q_F|Z_j]\}\big)^{1/q}+\|Z_j-\tilde{Z_j}\|_{L_q}\right)
      \\&\leq L\left(C_{\mathcal{N}}\|[\![ Z_j]\!]_r-Z_j\|_{L_q}+\|Z_j-\tilde{Z_j}\|_{L_q}\right) 
      \\&\leq LC_{\mathcal{N}}\|[\![ \tilde{Z_j}]\!]_r-Z_j\|_{L_q}+L\|Z_j-\tilde{Z_j}\|_{L_q}
      \\&\leq LC_{\mathcal{N}}\|[\![ \tilde{Z_j}]\!]_r-\tilde{Z_j}\|_{F}+L(C_\mathcal{N}+1)\|Z_j-\tilde{Z_j}\|_{L_q}.
      \end{align*}
      While the second term of the last inequality is bounded by~\eqref{eq: local_error_1}, we bound the first term by
      \begin{align*}
      \|[\![ \tilde{Z_j}]\!]_r-\tilde{Z_j}\|_{F}&\leq  \|[\![ \phi^{c_jh}_F(Y_i)]\!]_r-\phi_F^{c_jh}(Y_i)\|_F+\|\phi_F^{c_jh}(Y_i)-\tilde{Z_j}\|_F\\&\leq \|[\![ \phi^{c_jh}_F(Y_i)]\!]_r-\phi_F^{c_jh}(Y_i)\|_F+C_Lh^{\gamma_j+1}\\&\leq 
      (C_M\epsilon h+C_Lh^{\gamma_j+1}),
    \end{align*} 
    where we recall that the coefficient $c_j$ was used in the definition~\eqref{eq:stageorder} of stage order.
   In summary, we have
      \begin{align*}
       \|F(\mathcal{N}_j(Z_j))-F(\tilde{Z_j})\|_{L_q}&\leq  LC_{\mathcal{N}}(C_M\epsilon h+C_Lh^{\gamma_j+1})+L(C_{\mathcal{N}}+1)C_Zh(h\epsilon+h^{{\gamma}_2+1})
      \\&\leq  C_Fh(\epsilon +h^{\gamma_j}+h^{{\gamma}_2+1}).
    \end{align*} This concludes the proof of~\eqref{eq: F_diff} using the definition of $\tilde{\gamma}_j$.
\end{proof}
\begin{proof}[Proof of Theorem 3]
    By the triangular inequality, the local error satisfies
    \begin{equation} \label{eq:triangular}
        \|Y_{i+1}-\phi^h_{F}(Y_i)\|_{L_q}\leq \|Y_{i+1}-\tilde{Y}_{i+1}\|_{L_q}+\|\tilde{Y}_{i+1}-\phi^h_{F}(Y_i)\|_{L_q},
    \end{equation} where $\tilde{Y}_{i+1}=Y_i+h\sum^s_{j=1}b_j F(\tilde{Z}_j)$. Using that $\Omega_{s+1}$ and $\Psi_{s+1}$ are independent of $\Omega_1,\ldots \Omega_s$ and $\Psi_1,\ldots \Psi_s$, Theorem \ref{theorem: ny} yields
    \begin{align*}
       & \|Y_{i+1}-\tilde{Y}_{i+1}\|_{L_q} = \Big\|\mathcal{N}_{s+1}\Big(Y_{i}+h\sum^s_{j=1}b_jF(\mathcal N_j(Z_j))\Big)-\tilde{Y}_{i+1}\Big\|_{L_q}
       \\
         \leq & C_\mathcal{N}\Big\| \Big[\hspace{-.25cm}\Big[ Y_{i}+h\sum^s_{j=1}b_j F(\mathcal{N}_j(Z_j))\Big]\hspace{-.25cm}\Big]_r-Y_{i}-h\sum^s_{j=1}b_jF(\mathcal{N}_j(Z_j))\Big\|_{L_q}+\Big\|h\sum^s_{j=1}b_j(F(\mathcal{N}_j(Z_j))-F(\tilde{Z}_j))\Big\|_{L_q}
        \\
        \leq &C_\mathcal{N}\Big\|[\![ \tilde{Y}_{i+1}]\!]_r-Y_{i}-h\sum^s_{j=1}b_j F(\mathcal{N}_j(Z_j))\Big\|_{L_q}+\Big\|h\sum^s_{j=1}b_j(F(\mathcal{N}_j(Z_j))-F(\tilde{Z}_j))\Big\|_{L_q}
       \\
       \leq &C_{\mathcal{N}}\|[\![ \tilde{Y}_{i+1}]\!]_r-\tilde{Y}_{i+1}\|_{L_q}+(1+C_{\mathcal{N}})\Big\|h\sum^s_{j=1}b_j(F(\mathcal{N}_j(Z_j))-F(\tilde{Z}_j))\Big\|_{L_q}\\
        \leq & C\Big(\epsilon h+h^{\tau+1}+h\sum^s_{j=1}|b_j|\|(F(\mathcal{N}_j(Z_j))-F(\tilde{Z}_j))\|_{L_q}\Big).
    \end{align*} Hence, the local error is bounded by 
    \begin{equation*}
         \|Y_{i+1}-\phi^h_{F}(Y_i)\|_{L_q}\leq\|Y_{i+1}-\tilde{Y}_{i+1}\|_{L_q}+C_1h^{\tau+1}\leq Ch(\epsilon+h^{{\gamma}}+h^\tau).
    \end{equation*} This is turned into a bound on the global error using as in the proof of Theorem \ref{Thm: idealised lq}, which yields the bounds on the $L_q$ norm claimed in the statement of Theorem~\ref{Thm: stage order}. The tail bound is obtained using Markov's inequality; see Corollary~\ref{corollary: idealised_method_prob}.
\end{proof}

    \bibliographystyle{plain}
    \bibliography{Bibliography}

\end{document}